\documentclass[a4paper]{article}
\usepackage{amsmath, amssymb, amsthm, graphicx,
 bm, comment}
\usepackage[arrow,curve,matrix,arc,2cell]{xy}
\UseAllTwocells

\usepackage{oldgerm}
\usepackage{multicol}

\usepackage[pdftex,colorlinks,urlcolor=blue,pdfstartview=FitH]{hyperref}
\usepackage[left=3.70cm, right=3.70cm, top=3.70cm, bottom=3.70cm]{geometry}

\DeclareMathOperator{\Hom}{Hom}

\DeclareMathOperator{\im}{Im} 

\DeclareMathOperator{\Ker}{Ker}

\def\a{\alpha}

\newcommand{\ZZ}{\mathbb{Z}}
\newcommand{\mR}{\mathbb{R}}

\newcommand{\mP}{\ensuremath{\mathcal{P}}}

\newtheorem{thm}{Theorem}[section]
\newtheorem{lem}[thm]{Lemma}

\newtheorem{cor}[thm]{Corollary}
\newtheorem{prop}[thm]{Proposition}
\newtheorem{rmk}[thm]{Remark}

\theoremstyle{definition}
\newtheorem{defn}[thm]{Definition}

\def\lto{\longrightarrow}
\def\lmto{\longmapsto}
\def\leq{\leqslant}
\def\geq{\geqslant}
\def\Ai{$A_\infty$}

\newcommand{\m}{\mathfrak{m}}
\newcommand{\ba}{\ensuremath{\mathbf a}}
\newcommand{\bb}{\ensuremath{\mathbf b}}
\newcommand{\bc}{\ensuremath{\mathbf c}}

\date{}

\begin{document}

\title{Tensor product of filtered $A_\infty$-algebras}

\author{Lino Amorim}
\maketitle 

\begin{abstract}
We define the tensor product of filtered $A_\infty$-algebras, establish some of its properties and give a partial description of the space of bounding cochains in the tensor product. Furthermore we show that in the case of classical $A_\infty$-algebras our definition recovers the one given by Markl and Shnider. We also give a criterion that implies that a given $A_\infty$-algebra is quasi-isomorphic to the tensor product of two subalgebras. This will be used in a sequel to prove a K\"unneth Theorem for the Fukaya algebra of a product of Lagrangian submanifolds. 
\end{abstract}

\section{Introduction}

In this paper we study the tensor product of \Ai-algebras. The case of classical (or flat) \Ai-algebras, that is algebras with $\m_0=0$, was studied by several authors, see \cite{AmoTu}, \cite{Lod}, \cite{MarShn} and \cite{SanUmb, SanUmb2}. We will focus on a class of curved \Ai-algebras over the Novikov ring $\Lambda_0$, that has not been considered yet, namely filtered \Ai-algebras. These were introduced by Fukaya, Oh, Ohta and Ono in \cite{FOOO}, to study the obstruction to the existence of Lagrangian Floer cohomology and play an important role in symplectic geometry and homological mirror symmetry.

For classical $A_\infty$-algebras, it follows from the general theory of minimal models of operads (see \cite{operads}) that the tensor product of $A_\infty$-algebras exists. However this does not provide an explicit model for the tensor product. The first explicit construction of the tensor product of classical $A_\infty$-algebras was given by Saneblidze and Umble in \cite{SanUmb, SanUmb2}. They do it by constructing an explicit diagonal for the operad governing $A_\infty$-algebras, the associahedra. Later Markl and Shnider \cite{MarShn} gave a somewhat more conceptual construction of a diagonal which they claim coincides with the construction of Saneblidze-Umble. Markl and Shnider use a cubical decomposition  of the associahedra and then apply the usual Alexander-Whitney diagonal to this cubic complex. Loday \cite{Lod} gave another construction of a diagonal, this time using a simplicial decomposition of the associahedra. All these constructions are quasi-isomorphic but in fact it seems they give  the exact 
same $\m_k$ operations. This has not been checked for $k\geq 6$. 

We will take a different approach, first suggested in \cite{KonSoi}, that has the advantage of not using operads. This makes it easier to generalize to the filtered case.
First we reduce the problem to the case of filtered dg-algebras, that is filtered \Ai-algebras with $\m_k=0$ for $k\geq 3$. Given a (unital) filtered $A_\infty$-algebra $A$ we show there is a filtered dg-algebra $End_A$ quasi-isomorphic to $A$. For classical algebras this follows from the Yoneda embedding (see \cite{Sei}), but for filtered algebras it is new.
For filtered dg-algebras there is an obvious definition of tensor product, so we take 
$$A\otimes_\infty B:=End_A \otimes_{dg} End_B.$$
This definition is very natural  and satisfies (up to quasi-isomorphism) the usual properties of the tensor product of associative algebras. It might however seem a bit unsatisfactory as it does not give an $A_\infty$-algebra structure on the vector space $A\otimes B$. 
We can remedy this situation using the homological perturbation lemma to transfer the \Ai-algebra structure on $A\otimes_\infty B$ to the vector space $A\otimes B$. This way we obtain an explicit \Ai-algebra structure on $A \otimes B$, described in Theorem \ref{filteredtensor}. 
We summarize the main properties of the tensor product $\otimes_{\infty}$ in  
\begin{thm}\label{intprop}
 Let $A_i$ and $B_i$ be filtered $A_\infty$-algebras. We have the following:
\begin{enumerate}
\item If $A_1 \simeq A_2$ and $ B_1 \simeq B_2$ are quasi-isomorphic, then $A_1\otimes_\infty B_1 \simeq A_2 \otimes_\infty B_2$;
\item $A\otimes_\infty \mathbb{K} \simeq A$, where $\mathbb{K}$ is the ground field;
\item $A_1\otimes_\infty A_2\simeq A_2\otimes_\infty A_1$;
\item $A_1\otimes_\infty(A_2\otimes_\infty A_3)\simeq (A_1\otimes_\infty A_2)\otimes_\infty A_3$;
\item If $A_i$ are flat \Ai-algebras then $A_1\otimes_\infty A_2$ is quasi-isomorphic to the tensor product $A_1\otimes_{SU} A_2$ defined by Markl and Shnider in \cite{MarShn}. 
\end{enumerate} 
\end{thm}

Our second main theorem gives a set of conditions under which a filtered (or classical) \Ai-algebra $C$ is quasi-isomorphic to the tensor product of two subalgebras $A$ and $B$. For associative algebras the basic condition for this to happen is that elements of $A$ and $B$ should commute in $C$. Our definition of commuting subalgebras can be thought of as a generalization of this to the \Ai \ world. Although completely algebraic, this definition and the construction used to prove the following theorem arose in symplectic geometry, in the author's Ph.D. thesis \cite{Amothesis}. 
\begin{thm}\label{intcriterion}
Let $(A,\m^A)$ and $(B,\m^B)$ be commuting subalgebras of $(C,\mu)$ in the sense of Definitions \ref{subalgebra} and \ref{comsubalg}.
If $K:A\otimes B\lto C$ defined as $K(a\otimes b)=(-1)^{\vert a\vert}\mu_{2,0}(a,b)$ is an injective map which induces an isomorphism on $\mu_{1,0}$-cohomology then there is a (strict) quasi-isomorphism
$$A\otimes_\infty B \simeq C.$$ 
\end{thm}
The main application of this theorem, that we have in mind, is the proof of a K\"unneth Theorem for the Fukaya algebra.  Let $L$ be a Lagrangian submanifold of a symplectic manifold $(M,\omega)$, the Fukaya algebra $\mathcal{F}(L)$ is a filtered \Ai-algebra structure on the singular chain (or de Rham) complex of $L$, constructed in \cite{FOOO}. Given Lagrangian submanifolds $L_i\subset (M_i,\omega_i)$,  consider the product Lagrangian $L_1\times L_2 \subset (M_1\times M_2,\omega_1\oplus\omega_2)$. In \cite{Amo}, it is proved that $\mathcal{F}(L_1)$ and $\mathcal{F}(L_2)$ are commuting subalgebras of $\mathcal{F}(L_1 \times L_2)$. The above theorem then implies that $\mathcal{F}(L_1 \times L_2)$ is
quasi-isomorphic to $\mathcal{F}(L_1) \otimes_\infty \mathcal{F}(L_2)$.

For the applications in symplectic geometry it is important to understand the space of bounding cochains, or Maurer--Cartan elements, of a filtered \Ai-algebra. These are solutions of the equation
$$\sum_{k\geq 0}\m_k(x,\ldots, x)=\mP(x) e_A,$$
where $e_A$ is the unit of $A$ and $\mP(x)$ is some element in the Novikov ring $\Lambda_0$. We denote by $MC(A)$ the set of solutions to this equation, modulo an equivalence relation known as gauge equivalence. Given a bounding cochain $x$ we can deform $A$ to obtain an \Ai-algebra $(\hat{A},\m^{x})$ with $\m^{x}_0=\mP(x) e_A$. This is a classical \Ai-algebra over the Novikov field $\Lambda$. We prove the following

\begin{thm}\label{intbounding}
Let $A$ and $B$ be filtered \Ai-algebras. There is a map
 $$\boxtimes: MC(A) \times MC(B) \lto MC(A\otimes_{\infty}B),$$
 which satisfies $\mP(x \boxtimes y) = \mP(x)+\mP(y)$. When $A$ and $B$ are graded, connected \Ai-algebras, this map is a bijection. Moreover there is a (strict) quasi-isomorphism
 $$ (\hat{A},\m^x)\otimes_\infty (\hat{B},\m^y) \simeq (\widehat{A\otimes_\infty B}, \m^{\otimes,x \boxtimes y }).$$
\end{thm}

This paper is organized in the following way. In Section 2, we review the main aspects of the theory of filtered \Ai-algebras that we will use, mostly following \cite{FOOO}. In Section 3, we construct the filtered dg-algebra $End_A$, define the tensor product of filtered \Ai-algebras and establish its basic properties. In Section 4, we use the homological perturbation lemma to give a model for the tensor product $A \otimes_\infty B$ on the vector space $A\otimes B$. This will allow the comparison with the previous definitions of tensor product and prove the last part of Theorem \ref{intprop}. In Section 5, we define commuting subalgebras and prove Theorem \ref{intcriterion}. In the last section we study bounding cochains on the tensor product $A \otimes_\infty B$ and prove Theorem \ref{intbounding}.

\vspace{.3cm}

\noindent{\bf Acknowledgements:} This paper is a reinterpretation of some of the results in my Ph.D. thesis. I would like to thank my advisor Yong-Geun Oh for his continued help and support. I would also like to thank Dominic Joyce and Junwu Tu for useful comments.
During my Ph.D. I was partially supported by FCT through the scholarship  SFRH/ BD/30381/2006. During the preparation of this paper I was supported by EPSRC grant EP/J016950/1.

\vspace{.3cm}

\noindent{\bf Conventions:}
Given an homogeneous element $a$ in a graded module $A$, we will denote its degree by $\vert a \vert$. We will also use a shifted degree, $\vert\vert a \vert\vert = \vert a \vert -1$.

We use bold letters to denote elements in a tensor algebra,
$$\textbf{a}= a_1 \otimes \ldots \otimes a_n \in A^{\otimes n}$$
and $\vert\vert \textbf{a} \vert\vert = \sum^{n}_{i=1} \vert\vert a_i \vert\vert$.
Furthermore, we use Sweedler notation for the standard coproduct on a tensor coalgebra, namely:
$$ \Delta(\textbf{a})= \sum_{(\textbf{a})} {\bf a^{(1)}} \otimes {\bf a^{(2)}} = \sum_{i=0}^{n} (a_1\otimes \ldots \otimes a_i)\otimes (a_{i+1}\otimes \ldots \otimes a_n). $$

In this paper, a \emph{tree}  $T$ is always a planar tree, which is a finite graph with no cycles embedded in the plane, with a distinguished set of vertices with one incident edge. Of these one is chosen as the \emph{root} and the others are called \emph{leaves}. In all the figures these vertices are left implicit. We would like to point out that the tree might have other vertices with a single incident edge, other than the root and the leaves. All vertices but the root and the leaves are called \emph{internal} and we denote the set of these by $V(T)$. We denote by $E(T)$ the set of edges and by $E_{int}(T)$ the set of internal edges, that is not incident at the root or a leaf. Given a vertex $v$ there is a unique path on $T$ from $v$ to the root. The unique edge incident at $v$ that is part of this path is said to be \emph{outgoing} and the other edges incident at $v$ are \emph{incoming}. The \emph{valency} of $v$, $val(v)$ is the number of incoming edges.

\section{Filtered \Ai-algebras}
\subsection{\Ai-algebras and homomorphisms}
In this section we will review the definitions and main properties of filtered \Ai-algebras and homomorphisms introduced in \cite{FOOO}.
We start with the definition of a general \Ai-algebra.

\begin{defn}
 An \Ai-algebra over a ring $R$ consists of a $\ZZ_2$-graded $R$-module $A$ and a collection of multilinear maps $\m_{k}:A^{\otimes k} \lto A$ for each $k \geq 0$ of degree $k \pmod 2$ satisfying the following equation
\begin{align} \label{Ainf}
 \sum_{\substack{0\leq j\leq n\\1\leq i\leq n-j+1}}(-1)^{*}\m_{n-j+1}(a_1,\ldots,\m_{j}(a_{i},\ldots,a_{i+j-1}),\ldots,a_n)=0
\end{align}
where $*=\sum_{l=1}^{i-1} \vert\vert a_l \vert\vert$.

When $R$ is a field and $\m_0 =0$, we say $A$ is a classical \Ai-algebra.
The \Ai-algebra is graded if $A$ is a $\ZZ$-graded module and the maps $\m_k$ have degree $2-k$.
\end{defn}
\begin{rmk}
 Following the convention in \cite{Seisub}, from now on, given any expression which consists of a multilinear map applied to a block of entries, like (\ref{Ainf}) for example, we write $*$ for the sum of the shifted degrees of the entries lying to the left of the block.
\end{rmk}

In this paper we will be interested in a particular kind of \Ai-algebra defined over the Novikov ring. Let $\mathbb{K}$ be a field, the Novikov ring over $\mathbb{K}$ is defined as follows
$$\Lambda_0^{\mathbb{K}}=\Big\{\sum_{i=0}^\infty a_iT^{\lambda_i}\vert \lambda_i \in \mR, a_i\in \mathbb{K}, 0\leq\ldots\leq \lambda_i\leq\lambda_{i+1}\leq\ldots,\lim_{i\to\infty} \lambda_i= +\infty\Big\},$$
with maximal ideal
$$\Lambda_+^{\mathbb{K}}=\Big\{\sum_i a_iT^{\lambda_i}\vert \lambda_i > 0,\forall\ i\textrm{ with } a_i\neq0 \Big\}.$$ 
Localizing at $\Lambda_+^{\mathbb{K}}$ we obtain the Novikov field
$$\Lambda^{\mathbb{K}}=\Big\{\sum_{i=0}^\infty a_iT^{\lambda_i}\vert \lambda_i \in \mR, a_i\in\mathbb{K}, \lambda_i\leq\lambda_{i+1},\lim_{\lambda_i\to\infty}= +\infty\Big\}.$$
Note also that $\Lambda_0^{\mathbb{K}}$ has a natural filtration
$$F^\lambda\Lambda_0^{\mathbb{K}}=\Big\{\sum_i a_iT^{\lambda_i}\vert \lambda_i\geq \lambda,\forall\ i\textrm{ with } a_i\neq0 \Big\}.$$
Usually we will drop $\mathbb{K}$ from the notation and simply write $\Lambda_0$, $\Lambda$ and $\Lambda^+$.

Next, consider $G\subset \mathbb{R}_{\geq 0}\times2\mathbb{Z}$ and write $E:G\lto \mathbb{R}_{\geq 0}$ and $\mu:G\lto 2\mathbb{Z}$ for the natural projections. We say $G$ is a \emph{discrete submonoid}, if it is an additive submonoid satisfying
$$E^{-1}([0,c])\textrm{ is finite for any }c\geq 0.$$

\begin{defn}
 Let $G$ be a discrete submonoid, a $G$-gapped filtered \Ai-algebra $A=(A, \m)$ consists of a $\ZZ$-graded $\mathbb{K}$-vector space $A$ together with maps $\m_{k,\beta}:A^{\otimes k} \lto A$, for each $\beta \in G$ and $k\geq0$ of degree $2-k-\mu(\beta)$. These are required to satisfy $\m_{0,0}=0$ and for all $\beta \in G$ and homogeneous $a_1,\ldots,a_n \in A$:
\begin{align}\label{gapAinf}
\sum_{\substack{\beta_1+\beta_2=\beta\\0\leq j\leq n\\1\leq i\leq n-j+1}}(-1)^{*}\m_{n-j+1,\beta_2}(a_1,\ldots,\m_{j,\beta_1}(a_{i},\ldots,a_{i+j-1}),\ldots,a_n)=0.
\end{align}

We will say $(A,\m)$ is a filtered \Ai-algebra if it is $G$-gapped filtered for some $G$.
\end{defn}

Given a filtered \Ai-algebra $(A,\m)$, the tensor product $A\otimes_{\mathbb{K}}\Lambda_0$ inherits a filtration from $\Lambda_0$. We denote the completion of $A\otimes_{\mathbb{K}}\Lambda_0$ with respect to this filtration by $\hat{A}_0$ and define maps $\m_k:\hat{A}_0^{\otimes k} \lto \hat{A}_0$ by setting
$$  \m_k=\sum_{\beta\in G} \m_{k,\beta}T^{E(\beta)}.$$
The gapped condition ensures this is well defined and (\ref{gapAinf}) implies that $(\hat{A}_0, \m)$ is an \Ai-algebra over $\Lambda_0$.

Note that $(A,\m_{k,0})$ is a classical \Ai-algebra over $\mathbb{K}$. Therefore, by considering $G=\{ 0 \}$-gapped \Ai-algebras we recover the theory of classical \Ai-algebras over the field $\mathbb{K}$.

Next we introduce some terminology. 
\begin{defn}
 Let $(A,\m)$ be a $G$-gapped filtered \Ai-algebra.
\begin{itemize}
\setlength{\parsep}{0pt}
\setlength{\itemsep}{0pt}
\item[{\bf(a)}] If $\m_{k,\beta}=0$ when $k>2$ or when $k=2$ and $\beta\neq0$, we say $A$ is a filtered dg-algebra.
\item[{\bf(b)}] If $G \subset \mR_{\geq 0} \times \{0\}$, we say $A$ is graded.
\item[{\bf(c)}] If $\m_0=0$ we say $A$ is a flat \Ai-algebra.
\end{itemize}
\end{defn}

\begin{defn}
 An element $e_A\in A$ of degree $0$ is said to be a \emph{unit} of the filtered $A_\infty$-algebra $(A,\m)$ if
$$\m_{2,0}(e_A,a)=(-1)^{\vert a\vert}\m_{2,0}(a,e_A)=a$$ and $\m_{k,\beta}(\ldots,e,\ldots)=0$ for $(k,\beta)\neq (2,0)$.
\end{defn}

This notion of unit is sometimes called \emph{strict unit}. There are other, more flexible, notions of unit, namely cohomological unit or homotopy unit. As it is explained in \cite{Sei}, these notions are equivalent.

\begin{defn}
Let $(A,\m^A)$ and $(B,\m^B)$ be filtered $A_\infty$-algebras (for some discrete submonoid $G$). A \emph{filtered $A_\infty$-homomorphism} $F:A\lto B$ consists of a sequence of maps
$$F_{k,\beta}: A^{\otimes k}\lto B,\ k\geq 0,\ \beta\in G,$$
of degree $1-k-\mu(\beta)$, satisfying
\begin{itemize}
\setlength{\parsep}{0pt}
\setlength{\itemsep}{0pt}
\item[(a)] $F_{0,0}=0$,
\item[(b)] for each $n\geq 1$, $\beta\in G$ and $a_1,\ldots,a_n\in A$
\begin{align}\sum& \m_{l,\beta_0}^B(F_{k_1,\beta_1}(a_1,\ldots,a_{k_1}),\ldots,F_{k_l,\beta_l}(a_{n-k_l+1},\ldots,a_n))=\nonumber\\
&=\sum_{\substack{\beta_1+\beta_2=\beta\\0\leq j\leq n}}(-1)^{*}F_{n-j+1,\beta_2}(a_1,\ldots,\m^A_{j,\beta_1}(a_{i},\ldots,a_{i+j-1}),\ldots,a_n)\nonumber\end{align}
where the sum on the left hand side is taken over $l\geq 1,k_i\geq 0,\beta_i\in G$ such that $\beta_0+\beta_1+\ldots+\beta_l=\beta$ and $k_1+\ldots+k_l=n$.

For $n=0$ one must have 
$$\m^B_{0,\beta}+\sum_{l\geq 1} \m^B_{l,\beta_0}(F_{0,\beta_1},\ldots,F_{0,\beta_l})=\sum_{\beta_1+\beta_2=\beta}F_{1,\beta_1}(\m_{0,\beta_2}^A).$$
\end{itemize}
\end{defn}

\begin{defn}Let $F:A\lto B$ be a filtered \Ai-homomorphism.
\begin{itemize}
\item[{\bf(a)}] If $A$ and $B$ have units $e_A$ and $e_B$, we say $F$ is \emph{unital} if
$$F_{1,0}(e_A)=e_B\textrm{ and }F_{k,\beta}(\ldots,e_A,\ldots)=0\textrm{ for }(k,\beta)\neq (1,0).$$
\item[{\bf(b)}] We say $F$ is \emph{strict} if $F_{0,\beta}=0$ for all $\beta$.
\item[{\bf(c)}] We say $F$ is \emph{naive} if $F_{k,\beta}=0$ for all $(k,\beta)\neq (1,0)$.
\item[{\bf(d)}] If $F_{1,0}: A \lto B$ is an isomorphism we say $F$ is an \emph{isomorphism}. If the induced map on cohomology $ (F_{1,0})_*:H^*(A,\m_{1,0}^A)\lto H^*(B,\m_{1,0}^B)$ is an isomorphism we say $F$ is a \emph{quasi-isomorphism}.
\end{itemize}
\end{defn}

Given two filtered $A_\infty$-homomorphism $F: A\lto B$ and $G:B\lto C$ we can define  their composition as 
$$(G\circ F)_{n,\beta}(a_1,\ldots,a_n)=\sum G_{l,\beta_0}(F_{k_1,\beta_1}(a_1,\ldots,a_{k_l}),\ldots,F_{k_l,\beta_l}(\ldots,a_n))$$ where we sum over $l,\beta_i, k_i\geq 0$, such that $k_1+\ldots+k_l=n$ and $\beta_0+\beta_1+\ldots+\beta_l=\beta$.

In \cite{FOOO}, the authors develop the homotopy theory of filtered $A_\infty$-algebras, which we quickly review. Given a filtered $A_\infty$-algebra $(A,\m)$ one can define a filtered $A_\infty$-algebra $(A^{[0,1]},\m^{[0,1]})$ and naive filtered $A_\infty-$homomorphisms $i:A\lto A^{[0,1]}$, $p_{0}:A^{[0,1]}\lto A$ and $p_{1}: A^{[0,1]}\lto A$ satisfying:
\begin{itemize}
\item[(i)] $i$, $p_0$ and $p_1$ are quasi-isomorphism,
\item[(ii)] $p_1\circ i=p_0\circ i=id_A$,
\item[(iii)] $(p_1)_{1,0}\oplus (p_0)_{1,0}:A^{[0,1]}\lto A\oplus A$ is surjective.
\end{itemize}
These properties uniquely determine $A^{[0,1]}$ in an appropriate sense (see \cite[Chapter 4]{FOOO}). In Section 6 we will give a model for $A^{[0,1]}$ using the tensor product of filtered \Ai-algebras.

We can now give the following
\begin{defn}
Two filtered  $A_\infty$-homomorphisms $F_0,F_1:A\lto B$ are said to be \emph{homotopic} if there is a filtered $A_\infty$-homomorphism $F^{[0,1]}: A\lto B^{[0,1]}$ such that
$$F_0=p_0\circ F^{[0,1]}\textrm{ and }F_1=p_1\circ F^{[0,1]}.$$
A filtered $A_\infty$-homomorphism $F:A\lto B$ is called a \emph{homotopy equivalence} if there exists $G:B\lto A$ such that $F\circ G$ and $G\circ F$ are homotopic to the identity (the identity is the naive map defined by $(id)_{1,0}=id$).
\end{defn}

The following is one of the the most important theorems about $A_\infty$-algebras, it is sometimes referred to as the Whitehead theorem for $A_\infty$-algebras.
\begin{thm}[{\cite[Section 4.5]{FOOO}}]
Let $F:A\lto B$ be a filtered $A_\infty$-homomorphism. If $F$ is a quasi-homomorphism, it is a homotopy equivalence. If $F$ is unital (respectively strict), its homotopy inverse can also be taken to be unital (respectively strict).
\end{thm}

%\medskip

\subsection{Homological perturbation lemma} In this subsection we will review the homological perturbation lemma. This allows us to transfer the \Ai-algebra structure $(A,\m)$ to a different chain complex homotopic to $A$.

Consider the following situation: let $(A,\m)$ be a filtered $A_\infty$-algebra, $(V,d)$ be a chain complex and suppose we have chains maps $i:(V,d)\lto (A,\m_{1,0})$, $p:(A, \m_{1,0})\lto (V,d)$ and a homotopy $H: A\lto A$ satisfying,
\begin{align}\label{homotopy}
i\circ p-id_{A}=\m_{1,0}\circ H+H \circ \m_{1,0}.
\end{align}
We call such collection of maps \emph{homotopy data}. We have the following

\begin{thm}\label{hpl}
In the situation described, $V$ has the structure of filtered \Ai-algebra with $\mathfrak{\eta}_{1,0}=d$ and there is a filtered $A_\infty$-homomorphism $\varphi:(V, \eta)\lto(A,\m)$ with $\varphi_{1,0}=i$. Moreover there is a filtered \Ai -homomorphism $\psi: A \lto V$ with $\psi_{1,0}=p$ and a homotopy between $\varphi \circ \psi$ and $id_A$.
\end{thm}

A proof of this theorem for classical \Ai-algebras, together with explicit formulas for the maps is given in \cite{Mar}. For the case of filtered \Ai-algebras see \cite{FOOOcan}. Here we will only give a description of the maps $\eta_{k,\beta}$ and $\varphi_{k,\beta}$.

Let $\Gamma_k(G)$ denote the set of trees with $k$-leaves together with a map $\nu:V(T)\lto G$. We impose the condition that $\nu(v)\neq 0$ if $v$ is an internal vertex of valency zero or one.
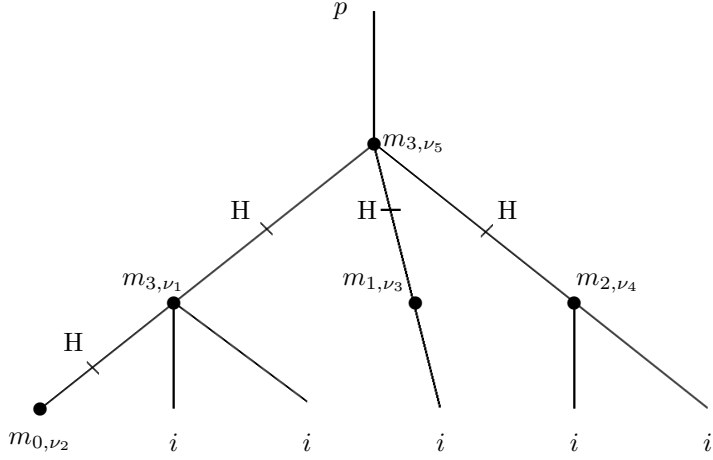
\begin{figure}
\begin{center}
 \setlength{\unitlength}{2.5pt}
\begin{picture}(100,70)(0,0)
\put(0,10){\circle*{2}}
\put(20,26){\circle*{2}}
\put(50,50){\circle*{2}}
\put(80,26){\circle*{2}}
\put(56.2,26){\circle*{2}}

\linethickness{0.3mm}
\put(50,50){\line(0,1){20}}
\linethickness{0.3mm}
\multiput(0,10)(0.15,0.12){333}{\line(1,0){0.15}}
\linethickness{0.3mm}
\put(20,10){\line(0,1){16}}
\linethickness{0.3mm}
\multiput(20,26)(0.16,-0.12){125}{\line(1,0){0.16}}
\linethickness{0.3mm}
\multiput(50,50)(0.12,-0.48){83}{\line(0,-1){0.48}}
\linethickness{0.3mm}
\multiput(50,50)(0.15,-0.12){333}{\line(1,0){0.15}}
\linethickness{0.3mm}
\put(80,10){\line(0,1){16}}
\linethickness{0.3mm}
\put(51,40){\line(1,0){3}}
\linethickness{0.3mm}
\multiput(33,38)(0.12,-0.12){15}{\line(1,0){0.12}}
\linethickness{0.3mm}
\multiput(7,17)(0.12,-0.12){15}{\line(1,0){0.12}}
\linethickness{0.3mm}
\multiput(66,36)(0.12,0.12){15}{\line(1,0){0.12}}

\put(45,70){\makebox(0,0)[cc]{$p$}}
\put(20,5){\makebox(0,0)[cc]{$i$}}
\put(40,5){\makebox(0,0)[cc]{$i$}}
\put(60,5){\makebox(0,0)[cc]{$i$}}
\put(80,5){\makebox(0,0)[cc]{$i$}}
\put(100,5){\makebox(0,0)[cc]{$i$}}
\put(70,40){\makebox(0,0)[cc]{H}}
\put(5,20){\makebox(0,0)[cc]{H}}
\put(30,40){\makebox(0,0)[cc]{H}}
\put(49,40){\makebox(0,0)[cc]{H}}
\put(85,29){\makebox(0,0)[cc]{$m_{2,\nu_4}$}}
\put(56,50){\makebox(0,0)[cc]{$m_{3,\nu_5}$}}
\put(50,29){\makebox(0,0)[cc]{$m_{1,\nu_3}$}}
\put(0,5){\makebox(0,0)[cc]{$m_{0,\nu_2}$}}
\put(17,29){\makebox(0,0)[cc]{$m_{3,\nu_1}$}}
\end{picture}
\end{center}
\caption{Example of an element of $\Gamma_5(\sum_{i=1}^{5}\nu_i)$. Note that $\nu_3,\ \nu_2\neq 0$.}
\label{figexgamma}
\end{figure}

Now, given $\beta\in G$, consider the set $\Gamma_k(\beta)=\{T \in \Gamma_k(G) \vert \sum_{v \in V(T)}\nu(v)=\beta\}$. We can easily see that $\Gamma_k(\beta)$ is finite. For each $T\in \Gamma_k(\beta)$, we obtain $\overline T$ from $T$ by inserting one additional vertex into each internal edge of $T$. To each inserted vertex assign the map $H$. We use $T$ as a flow chart to define a map
\begin{align}\label{hpformula}
 \eta_T: V^{\otimes k}\lto V.
\end{align}
We assign to each $v\in V(T)$ the map $\m_{val(v),\nu(v)}$ and we assign $p$ to the root and $i$ to the leaves. For example, the tree in Figure \ref{figexgamma} gives the map
\begin{align}
\eta_T(&v_1,v_2,v_3,v_4,v_5)=\nonumber\\
&=p\big(\m_{3,\nu_5}(H\circ \m_{3,\nu_1}(H(\m_{0,\nu_2}),i(v_1),i(v_2)),H\circ \m_{1,\nu_3}(i(v_3)),H\circ \m_{2,\nu_4}(i(v_4),i(v_5)))\big).\nonumber\end{align}
Then we define
$$\eta_{k,\beta}=\sum_{T\in\Gamma_k(\beta)}\eta_T.$$

The construction of the map $\varphi: V\lto A$ is similar. We take $\varphi_{k,\beta}=\sum_{T\in\Gamma_k(\beta)}\varphi_T$. The map $\varphi_T$ is defined in the same way as $\eta_T$, the only difference is that we assign $H$ to the root vertex (instead of $p$ as in the case of $\eta_T$).

When $A$ is unital, we can ensure $V$ and $\varphi$ are also unital by imposing side conditions on $p$, $i$ and $H$.

\begin{prop}\label{sideconditions}
Suppose $(A,\m)$ is unital and that the homotopy data satisfies
$$p\circ i = id, \ H\circ i=0, \  p\circ H=0 \ \textrm{and} \ H^2=0. $$
Additionally, assume there is an element $e_V$ such that $i(e_V)=e_A$. Then $e_V$ is a unit for $(V,\eta)$ and $\varphi$ is a unital quasi-isomorphism.
\end{prop}
\begin{proof} Simply inspecting the formulas given above for $\eta_{k,\beta}$ and applying the side conditions we see that 
$$\eta_{k+1,\beta}(v_1,\ldots,e_V,\ldots,v_k)= p(\m_{k+1,\beta}(i(v_1),\ldots,e_A,\ldots,i(v_k))),$$
which readily implies $e_V$ is a unit. The same argument shows $\varphi$ is unital.
The first condition $p\circ i = id$ together with the homotopy imply that $i=\varphi_{1,0}$ is an isomorphism in cohomology. Thus, we conclude $\varphi$ is a unital quasi-isomorphism.
\end{proof}

The standard application of the homological perturbation lemma is the following proposition proved in \cite{FOOOcan}.

\begin{prop}
 Any filtered \Ai-algebra $(A, \m)$ is quasi-isomorphic to a filtered \Ai-algebra $(A', \m')$ with $\m'_{1,0}=0$. We call $(A', \m')$ a canonical model for $A$.
\end{prop}
\begin{proof}
We pick a subspace $W\subseteq A$ such that $W \oplus \Ker \m_{1,0} = A$. Since $\m_{1,0}^2=0$ we can choose $H \subseteq \Ker \m_{1,0}$ satisfying $H \oplus \m_{1,0}(W)= \Ker \m_{1,0} $. This gives a decomposition 
$$A=H \oplus W \oplus \m_{1,0}(W).$$
Note that if $H$ is nontrivial we can choose it so that $e_A \in H$. Using the decomposition we define a inclusion $i:H \lto A$, a projection $p:A \lto H$ and a map $H:A\lto A$ of degree $-1$, which is zero when restricted to $H\oplus W$ and is the inverse of $\m_{1,0}$ on $\m_{1,0}(W)$. It is straightforward to check that this defines homotopy data satisfying the side conditions. Applying the homological perturbation lemma to this data we obtain a filtered \Ai-algebra $(H, \eta_{k,\beta})$ with $\eta_{1,0}=0$ quasi-isomorphic to $A$.
\end{proof}

\subsection{Bounding cochains}

Here we will review some basic definitions and facts about bounding cochains. For details see \cite[Section 4.3]{FOOO}.

\begin{defn} Let $(A,\m)$ be a ($G$-gapped) filtered $A_\infty$-algebra and consider $x=\sum_{i=0}^{\infty} x_i T^{\lambda_i} \in \hat{A}_0$ with $\lambda_i \geq \lambda$ for some $\lambda>0$ and $x_i$ of odd degree. We say $x$ is a \emph{bounding cochain} if there is $\mP(x)\in\Lambda_0$ such that
\begin{align}\label{maurer_cartan}
 \sum_{k\geq 0}\m_k(x,\ldots, x)=\mP(x) e_A.
\end{align}
\end{defn}
%We call $\mP$ the \emph{potential}. 
Note that the authors in \cite{FOOO} use the term \emph{ weak bounding cochain} and reserve the term \emph{bounding cochain} to the case when $\mP(x)=0$. Bounding cochains are also called Maurer--Cartan elements. We denote by $\widehat{MC}(A)$ the set of all bounding cochains. 

When $A$ is graded, we require $|x_i|=1$ for all $i$. In this case $\mP(x)=0$ since the left hand side of (\ref{maurer_cartan}) has degree $2$ and $e_A$ has degree $0$.

Given a bounding cochain $x$ we deform the $A_\infty$-algebra $A$ in the following way. Consider $\gamma_i=(\lambda_i,1-\vert x_i \vert) \in \mathbb{R}_{\geq 0}\times2\mathbb{Z}$ and let $\bar{G}$ be the smallest submonoid containing $G$ and all the $\gamma_i$. This is a discrete submonoid. We then define a new $\bar{G}$-gapped $A_\infty$-algebra $(A,\m^x)$, by setting
$$\m_k^x(a_1,\ldots,a_k)=\sum_{i_0,\ldots,i_k}\m_{k+i_0+\ldots+i_k}(x,\ldots,x,a_1,x,\ldots,x,a_k,x,\ldots,x).$$
Observe that $(A,\m^x)$ is graded whenever $A$ is graded. 

Since $\m_0^x=\mP(x)e_A$ the $A_\infty$-relation implies $\m_1^x : \hat{A}_0 \lto \hat{A}_0$ has square zero. Thus given a bounding cochain $x$ we can define the cohomology for the pair $(A,x)$:

$$H^*(A,x;\Lambda_0):=\frac{\Ker \m_1^x}{\im \ \m_1^x}.$$

Since it is often more convenient to work over a field, we consider the \Ai-algebra $(\hat{A}, \m^x)$ where $\hat{A}=\hat{A}_0 \otimes_{\Lambda_0}\Lambda$. Observe that this is a classical \Ai-algebra over the field $\Lambda$, even when $\mP(x) \neq 0$ (since $\m_0^x=\mP(b)e_A$, terms in the \Ai-equations involving $\m_0^x$ cancel).
Hence we can think of $(\hat{A}, \m^x)$ as a classical \Ai-algebra with one additional invariant $\mP(x) \in \Lambda_0$.

Next we introduce an equivalence relation on $\widehat{MC}(A)$. We say $x_0,x_1\in\widehat{MC}(A)$ are \emph{gauge equivalent} if  there is $\bar x \in \widehat{MC}(A^{[0,1]})$ such that $p_0(\bar{x})=x_0$ and $p_1(\bar{x})=x_1$, in which case we write $x_1 \sim x_0$. 
We denote the set of equivalence classes by $MC(A)$. In fact, if $x_1\sim x_0$, then $\mP(x_1)=\mP(x_0)$ so that $\mP$ is a well defined function on $MC(A)$ (see \cite[Section 4.3]{FOOO}). Moreover $H^*(A,x_0;\Lambda_0) \cong H^*(A, x_1;\Lambda_0)$.

The following theorem is proved in Sections 4.3 and 5.2 of \cite{FOOO}:

\begin{thm} If $F:(A,\m^A)\lto(B,\m^B)$ is a unital quasi-isomorphism, the induced map $F_*: MC(A)\lto MC(B)$ is a bijection. Moreover 
$\mP(F_*(x))=\mP(x)$,
$$H^*(A,x;\Lambda_0) \cong H^*(B,F_*(x);\Lambda_0)$$
and $(\hat{A}, \m^x)$ is quasi-isomorphic to $(\hat{B}, \m^{F_*(x)})$, as classical \Ai-algebras over $\Lambda$.
\end{thm}

\section{Tensor product of filtered \Ai-algebras}

\subsection{Definition of the tensor product}
We will now consider the problem of defining the tensor product of filtered $A_\infty$-algebras. We will first prove that any \Ai-algebra is quasi-isomorphic to a filtered dg-algebra and then define the tensor product for these. This approach is inspired by the definition of tensor product of classical \Ai-algebras in \cite{KonSoi}.

Let $(A, \m)$ be a unital filtered \Ai-algebra. We define $End_A$ to be the subspace of $$\Hom_{\mathbb{K}}\big(\bigoplus_{s\geq 0} A\otimes A^{\otimes s}, A\big)$$ of elements satisfying $\rho_s(\bullet,\ldots, e_A,\ldots)=0$. An element $\rho = \{\rho_s\}_s \in End_A$ has degree $\vert\rho\vert$ if each $\rho_s$ has degree $\vert \rho \vert -s$. For each $\beta \in G$ we define the operations
\begin{align}
(\mu_{0,\beta})_s(v,a_1,\ldots,a_s)&=\sum_{\beta_1+\beta_2=\beta}\m_{s+2,\beta_1}(\m_{0,\beta_2},v,a_1,\ldots,a_s),\nonumber\end{align}

\begin{align}
(\mu_{1,\beta}(\rho))_s(v,a_1,\ldots,a_s)&=\sum_{i\geq 0}-\m_{s-i+1,\beta}(\rho_i(v,a_1,\ldots,a_i),\ldots,a_s)\nonumber\\
&+\sum_{i\geq 0}(-1)^{\vert\rho\vert}\rho_{s-i}(\m_{i+1,\beta}(v,a_1,\ldots,a_i),\ldots,a_s)\nonumber\\
&+\sum_{j,i\geq 0}(-1)^{\vert\rho\vert+*}\rho_{s-j+1}(v,\ldots,\m_{j,\beta}(a_{i+1},\ldots,a_{i+j}),\ldots,a_s),\nonumber\\
(\mu_{2,0}(\rho,\tau))&=(-1)^{\vert\rho\vert} \rho \circ \tau \nonumber\\
\textrm{where} \ 
(\rho \circ \tau)_s(v,a_1,\ldots,a_s)&=\sum_{i\geq 0}\rho_{s-i}(\tau_i(v,a_1,\ldots,a_i),\ldots,a_s).\nonumber
\end{align}
For all the other $(k,\beta)$ we define $\mu_{k,\beta}=0$. 

When $A$ is a flat \Ai-algebra, one can regard $A$ as a (right) \Ai-module over itself. In this case, $End_A$ is simply the filtered \Ai-algebra of \Ai-pre-homomorphisms of this module (see \cite{Seisub}). In general, $A$ is not a right module over itself, but as we will see, $End_A$ is still a filtered \Ai-algebra.

\begin{lem}\label{EndAdga}
$End_A=(End_A,\mu)$ is an unital filtered dg-algebra.
\end{lem}
\begin{proof}
The proof consists of a series of long and tedious computations. We start with
\begin{align}
&\left[\sum_{\beta_1+\beta_2=\beta}\mu_{1,\beta_1}(\mu_{0,\beta_2})\right]_s(v,a_1,\ldots,a_s)=\nonumber\\
&=\sum_{\beta_1+\beta'+\beta''=\beta}-\m_{s-i+1,\beta_1}(m_{i+2,\beta'}(\m_{0,\beta''}, v,\ldots,a_i),\ldots,a_s)\nonumber\\
&+\sum_{\substack{i\geq 0\\\beta_1+\beta'+\beta''=\beta}}\m_{s-i+2,\beta'}(\m_{0,\beta''},\m_{i+1,\beta_1}(v,\ldots,a_i),\ldots,a_s)\nonumber\\
&+\sum_{\substack{i,j\geq 0\\\beta_1+\beta'+\beta''=\beta}}(-1)^{* +1} \m_{s+j+1,\beta'}(\m_{0,\beta''},v,\ldots,\m_{j,\beta_1}(a_{i+1},\ldots,a_{i+j}),\ldots,a_s).\nonumber
\end{align}
By the $A_\infty$-equation this equals
\begin{align}
\sum_{\beta''}\sum_{\beta_1+\beta'=\beta-\beta''}\m_{s+3,\beta'}(\m_{0,\beta_1},\m_{0,\beta''},v,\ldots,a_s)-\m_{s+3,\beta'}(\m_{0,\beta''},\m_{0,\beta_1},v,\ldots,a_s)=0.\nonumber
\end{align}So we conclude $\sum_{\beta_1+\beta_2=\beta}\mu_{1,\beta_1}(\mu_{0,\beta_2})=0$.
The next equation can be handled in the same way:
\begin{align}
\sum_{\beta_1+\beta_2=\beta}\mu_{1,\beta_1}&(\mu_{1,\beta_2}(\rho))_s(v,a_1,\ldots,a_s)=\sum_{i\geq 0}-\m_{s-i+1,\beta_1}((\mu_{1,\beta_2}(\rho))_i(v,\ldots,a_i),\ldots,a_s)+\nonumber\\
&+\sum_{i\geq 0}(-1)^{\vert\rho\vert+1}(\mu_{1,\beta_2}(\rho))_{s-i}(\m_{i+1,\beta_1}(v,\ldots,a_i),\dots,a_s)+\nonumber\\
&+\sum_{i,j\geq 0}(-1)^{\vert\rho\vert+*}(\mu_{1,\beta_2}(\rho))_s(v,\ldots,\m_{j,\beta_1}(a_{i+1},\ldots,a_{i+j}),\ldots,a_s).\nonumber
\end{align}
Further expanding and canceling all the terms we see this is equal to
\begin{align}
&\sum_{j\geq0}\sum_{\substack{\beta_1+\beta_2=\beta\\i\geq j}}\m_{s-i+1,\beta_1}(\m_{i-j+1,\beta_2}(\rho_j(v,\ldots,a_j),\ldots,a_i)\ldots a_s)\nonumber\\
&+\sum_{j\geq 0}\sum_{\substack{\beta_1+\beta_2=\beta\\l\geq j\\i\geq 0}}(-1)^{\epsilon}\m_{s-i-j+2,\beta_1}(\rho_j(v,\dots,a_j),\ldots,\m_{i,\beta_2}(a_{l+1},\ldots,a_{l+i}),\ldots,a_s)\nonumber\\
&-\sum_j\sum_{\substack{\beta_1+\beta_2=\beta\\0\leq i\leq j}}\rho_{s-j}(\m_{j-i+1,\beta_1}(\m_{i+1,\beta_2}(v,\ldots,a_i)\ldots,a_{j}),\ldots,a_s)\nonumber\\
&-\sum_j\sum_{\substack{\beta_1+\beta_2=\beta\\l\geq 0\\i\geq 0}}(-1)^{*}\rho_{s-j}(\m_{j-i+2,\beta_1}(v,\ldots,\m_{i,\beta_2}(a_{l+1},\ldots,a_{l+i}),\ldots,a_{j}),\ldots,a_s),\nonumber
\end{align}
where $\epsilon=\vert\rho\vert+\vert v \vert +\sum_{p=1}^l \vert\vert a_p\vert\vert$. Adding the first two terms and the last two, and using the $A_\infty$-equation, we obtain
\begin{align}
-\sum_{\substack{\beta_1+\beta_2=\beta}}\sum_j&\m_{s-j+2,\beta_1}(\m_{0,\beta_2},\rho(v,\ldots,a_j),\ldots,a_s)+\nonumber\\
&+\sum_{\substack{\beta_1+\beta_2=\beta}}\sum_j\rho_{s-j}(\m_{j+2,\beta_1}(\m_{0,\beta_2},v,\dots,a_j),\ldots,a_s)=\nonumber
\end{align}
\begin{align}
=\left[-\mu_{2,0}(\mu_{0,\beta},\rho)+(-1)^{\vert\rho\vert}\mu_{2,0}(\rho,\mu_{0,\beta})\right]_s(v,a_1,\ldots,a_s).\nonumber
\end{align}
So we conclude that
$$\sum_{\substack{\beta_1+\beta_2=\beta}}\mu_{1,\beta_1}(\mu_{1,\beta_2}(\rho))+\mu_{2,0}(\mu_{0,\beta},\rho)+(-1)^{\vert\vert\rho\vert\vert}\mu_{2,0}(\rho,\mu_{0,\beta})=0.$$
The last two equations we have to check
$$\mu_{2,0}(\mu_{1,\beta}(\rho),\tau)+(-1)^{\vert\vert\rho\vert\vert}\mu_{2,0}(\rho,\mu_{1,\beta}(\tau))+\mu_{1,\beta}(\mu_{2,0}(\rho,\tau))=0,$$
and
$$\mu_{2,0}(\mu_{2,0}(\rho,\tau),\theta)+(-1)^{\vert\vert\rho\vert\vert}\mu_{2,0}(\rho,\mu_{2,0}(\tau,\theta))=0,$$
simply state that (up to a change in sign) $\mu_{2,0}$ is associative and $\mu_{1,\beta}$ is a derivation of $\mu_{2,0}$. Since these do not involve $\mu_{0,\beta}$, they can be derived from the flat case (see \cite{Seisub}).

Finally we can easily verify that
\begin{align}
(id)_r(v,a_1,\ldots,a_r)=\left\{\begin{array}{ll}
v& ,r=0\\
0& ,\textrm{otherwise}
\end{array}\right.\nonumber
\end{align} defines a unit in $End_A$.
\end{proof}

Next we will show that $End_A$ is quasi-isomorphic to $A$. We begin with the following proposition:

\begin{prop}
Let $A$ be a unital filtered $A_\infty$-algebra and $End_A$ be the filtered dg-algebra described above. Consider $F_{k,\beta}:A^{\otimes k}\lto End_A$ defined as 
$$F_{k,\beta}(\a_1,\ldots,\a_k)_s(v,a_1,\ldots,a_s)=\m_{k+s+1,\beta}(\a_1,\ldots,\a_k,v,a_1,\ldots,a_s),$$ for $k>0$. This defines a strict $A_\infty$-homomorphism.
\end{prop}
\begin{proof}
First note that
$$\mu_{0,\beta}=\sum_{\beta_1+\beta_2=\beta}F_{1,\beta_1}(\m_{0,\beta_2})$$
by definition. Next we need to check that
\begin{align}
\sum_{\substack{\beta_1+\beta_2=\beta\\1\leq i\leq k-1}}&(-1)^{ \vert F_i(\a_1,\ldots,\a_i)\vert}F_{i,\beta_1}(\a_1,\ldots,\a_{i}) \circ F_{k-i,\beta_2}(\a_{i+1},\ldots,\a_{k}) \nonumber\\
&+\mu_{1,\beta_1}\big(F_{k,\beta_2}(\a_1,\ldots,\a_{k})\big)\nonumber\\
&=\sum_{\substack{\beta_1+\beta_2=\beta\\0\leq j\leq k}} (-1)^*F_{k-j+1,\beta_1}\big(\a_1,\ldots,\m_{j,\beta_2}(\a_{i+1},\ldots,\a_{i+j})\ldots,\a_k\big).\nonumber
\end{align}
If we denote by $L$ and $R$ the left and right-hand side of the above equation, we need to check $L_s(v,a_1,\ldots,a_s)=R_s(v,a_1,\ldots,a_s)$.
Noting that $\vert F_i(\a_1,\ldots,\a_i)\vert=\sum_{l=1}^{i}\vert\vert\a_l\vert\vert +1$, we compute
\begin{align}
L_s&(v,a_1,\ldots,a_s)=\sum_{0\leq j\leq s}-\m_{s-j+1,\beta_1}\big(\m_{k+j+1,\beta_2}(\a_1,\ldots,\a_k,v,a_1,\ldots,a_j),\ldots,a_s\big)\nonumber\\
&+\sum_{0\leq j\leq s}(-1)^{*+1}\m_{k+s-j,\beta_2}\big(\a_1,\ldots,\a_k,\m_{j+1,\beta_1}(v,a_1,\ldots,a_j),\ldots,a_s\big)\nonumber\\
&+\sum_{\substack{0\leq j\leq s\\ 0\leq l\leq s-j}}(-1)^{*+1}\m_{k+s-j,\beta_2}\big(\a_1,\ldots,\a_k,v,\ldots,\m_{j,\beta_1}(a_{l+1},\ldots,a_{l+j}),\ldots,a_s\big)\nonumber\\
&+\sum_{\substack{0\leq j\leq s\\ 1\leq i\leq k-1}}(-1)^{*+1}\m_{i+s-j+1,\beta_1}\big(\a_1,\ldots,\a_i,\m_{k-i+j+1,\beta_2}(\a_{i+1},\ldots,v,\ldots,a_{j}),\ldots,a_s\big)\nonumber
\end{align}
and
$$R_s(v,a_1,\ldots,a_s)=\sum_{\substack{0\leq j\leq k\\ 0\leq i\leq k-j}}(-1)^{*}\m_{s+k-j,\beta_1}\big(\a_1,\ldots\m_{j,\beta_2}(\a_{i+1},\ldots,\a_{i+j}),\ldots,v,\ldots,a_s\big).$$
Thus we conclude that $L_s(v,a_1,\ldots,a_s)-R_s(v,a_1,\ldots,a_s)=0$, since this is simply the $A_\infty$-equation on $A$.
\end{proof}

Next we observe that $F_{1,0}$ induces an isomorphism on cohomology by defining a homotopy inverse. Consider the map $P:End_A\lto A$ defined by $P(\rho)=(-1)^{\vert\rho\vert}\rho_0(e_A)$. We can easily see that $P:(End_A,\mu_{1,0})\lto(A,\m_{1,0})$ is a chain map and $P\circ F_{1,0}=id_A$. 

\begin{lem}\label{homotopytensor}
Define $H:End_A \lto End_A$ by
$$(H\rho)_s(v,a_1,\ldots,a_s)=(-1)^{\vert\rho\vert}\rho_{s+1}(e_A,v,a_1,\ldots,a_s).$$
Then
$F_{1,0}\circ P - id_{End_A}=\mu_{1,0}\circ H+H\circ\mu_{1,0}$ and $F_{1,0}$ is a quasi-isomorphism.
\end{lem}
\begin{proof}
 We compute
 \begin{align}(\mu_{1,0} \circ H(\rho))_s(v,a_1,\ldots&,a_s)=\sum_{0\leq i\leq s}(-1)^{\vert\rho\vert+1}\m_{s-i+1,0}(\rho_{i+1}(e,v,a_1,\ldots,a_i),\ldots,a_s)\nonumber\\
 +&\sum_{0\leq i\leq s}\rho_{s-i}(e,\m_{i+1,0}(v,a_1,\ldots,a_i),\ldots,a_s)\nonumber\\
 +&\sum_{\substack{1\leq j\leq s\\0\leq i\leq s-j}}(-1)^{* +1}\rho_{s-j}(e,v,a_1,\ldots,\m_{j,0}(a_{i+1},\ldots,a_{i+j}),\ldots,a_s).\nonumber
 \end{align}
 Similarly
  \begin{align}(H\circ \mu_{1,0}(\rho))_s(v,a_1,\ldots&,a_s)=\sum_{0\leq i\leq s}(-1)^{\vert\rho\vert}\m_{s-i+1,0}(\rho_{i+1}(e,v,a_1,\ldots,a_i),\ldots,a_s)\nonumber\\
 +&(-1)^{\vert\rho\vert}\m_{s+2,0}(\rho_0(e),v,\ldots,a_s)\nonumber\\
 +&\sum_{0\leq i\leq s}-\rho_{s-i}(\m_{i+2,0}(e,v,a_1,\ldots,a_i),\ldots,a_s)\nonumber\\
 +&\sum_{\substack{1\leq j\leq s\\0\leq i\leq s-j}}(-1)^{*}\rho_{s-j}(e,v,\ldots,\m_{j,0}(a_{i+1},\ldots,a_{i+j}),\ldots,a_s)\nonumber\\
 +&\sum_{0\leq j\leq s}-\rho_{s-j+1}(e,\m_{j+1,0}(v,a_1,\ldots,a_j),\ldots,a_s).\nonumber
 \end{align}
 Comparing the two expressions we obtain
 \begin{align}
 \big((\mu_{1,0}H&+H\mu_{1,0})\rho\big)_s(v,a_1,\ldots,a_s)=\nonumber\\
 &=(-1)^{\vert\rho\vert}\m_{s+2,0}(\rho_0(e),v,\ldots,a_s)-\sum_{i\geq 0}\rho_{s-i}(\m_{i+2,0}(e,v,\ldots,a_i),\ldots,a_s)\nonumber\\
 &=(-1)^{\vert\rho\vert}\m_{s+2,0}(\rho_0(e),v,\ldots,a_s)-\rho_s(v,a_1,\ldots,a_s)\nonumber\\
 &=\Big[\big((F_{1,0}\circ P)\rho\big)-\rho\Big]_s(v,a_1,\ldots,a_s).
 \end{align}
 The second equality holds because $A$ is unital, so $\m_{i+2,0}(e,v,a_1,\ldots,a_i)=0$ unless $i=0$.

Since $P$ is a homotopy inverse of $F_{1,0}$, the conclusion follows.
\end{proof}

\begin{cor}\label{EndA}
Let $A$ be a unital filtered $A_\infty$-algebra. Then $A$ is quasi-isomorphic to the filtered dg-algebra $End_A$.
\end{cor}

We are now ready to define the tensor product of filtered \Ai-algebras. Given $A$ and $B$ filtered $G_A$ and $G_B$-gapped $A_\infty$-algebras, we will define their tensor product as a $G$-gapped $A_\infty$-algebra, for $G=G_A+G_B$. In fact the tensor product is bi-gapped in the following sense: for each $\beta_1\in G_A$ and $\beta_2\in G_B$ we will define $\m_{k,\beta_1\times\beta_2}^\otimes$ and given $\beta\in G_A+G_B$ we take
$$\m_{k,\beta}^\otimes=\sum_{\substack{\beta_1\times\beta_2\in G_A\times G_B\\\beta_1+\beta_2=\beta}}\m_{k,\beta_1\times\beta_2}^\otimes.$$

We start by defining the tensor product of filtered dg-algebras.

\begin{lem}\label{tensordg}
 Let $(C,\mu^{C})$ and $(D,\mu^{D})$ be unital filtered dg-algebras. On the vector space $C\otimes_{\mathbb{K}}  D$ consider the operations
\begin{align}
\mu^{\otimes}_{0,\beta_1\times\beta_2}&=\left\{\begin{array}{ll}
\mu_{0,\beta_1}^C\otimes e_D& \beta_2=0\\
e_C\otimes\mu_{0,\beta_2}^D& \beta_1=0\\
0&  \beta_1, \beta_2 \neq 0
\end{array},\right.\nonumber\end{align}
\begin{align}
\mu^{\otimes}_{1,\beta_1\times\beta_2}(c\otimes d)&=\left\{\begin{array}{ll}
\mu_{1,\beta_1}^C(c)\otimes d& \beta_2=0\\
(-1)^{\vert c\vert}c\otimes\mu_{1,\beta_2}^D(d)& \beta_1=0\\
0& \beta_1, \beta_2 \neq 0
\end{array},\right.\nonumber\\
\mu^{\otimes}_{2,0}(c_1\otimes d_1,c_2\otimes d_2)&=(-1)^{\vert c_2\vert\vert d_1\vert}\mu_{2,0}^C(c_1,c_2)\otimes \mu_{2,0}^D(d_1,d_2).\nonumber
\end{align}
Then $(C\otimes D, \mu^{\otimes})$ is an unital filtered dg-algebra, which we denote by $C\otimes_{dg} D$
\end{lem}
\begin{proof}
We omit the simple verification of the \Ai equations. It is also clear that $e_A\otimes e_B$ is a unit.
\end{proof}

\begin{rmk}
 To avoid confusion, we spell out the notation above when $\beta=0$: 
$$\mu^{\otimes}_{1,0}(c\otimes d)= \mu_{1,0}^C(c)\otimes d + (-1)^{\vert c\vert}c\otimes\mu_{1,0}^D(d).$$
\end{rmk}

We are now ready to state our main definition.
\begin{defn} Let $A$ and $B$ be filtered $A_\infty$-algebras. We define their tensor product as 
$$A\otimes_\infty B:=End_A\otimes_{dg} End_B.$$
\end{defn}

If $A$ and $B$ are classical \Ai-algebras over $\mathbb{K}$, then they are also $G=\{0\}$-gapped \Ai-algebras, then the above also gives a definition of classical \Ai-algebras. For these algebras the tensor product was already defined, see \cite{MarShn, SanUmb, Lod}. All these constructions define the tensor product as an $A_\infty$-algebra structure on the vector space $A\otimes B$. Our definition uses a different underlying vector space and therefore is not directly comparable. In the next section, we will remedy this by using the homological perturbation lemma to transfer our $A_\infty$-structure on $A\otimes_\infty B$ to the vector space $A\otimes B$. 
%This will also allow us to compare our definition with the ones already in the literature.

\subsection{Some properties of the tensor product}\label{secproperties}

Let $A_1$, $A_2$, $B_1$ and $B_2$ be filtered $A_\infty$-algebras and let $\varphi: A_1\lto A_2$ and $\psi: B_1\lto B_2$ be filtered $A_\infty$-homomorphisms. We will show that these induce a filtered $A_\infty$-homomorphism
$$\varphi\otimes_\infty\psi:A_1\otimes_\infty B_1\lto A_2\otimes_\infty B_2.$$
Although not functorial, this construction satisfies the additional property that if $\varphi$ and $\psi$ are quasi-isomorphisms then $\varphi\otimes_\infty\psi$ is also a quasi-isomorphism. This will be enough to show that the tensor product is well defined and satisfies some monoidal properties.

Following our strategy, we begin with the case of filtered dg-algebras. We will use the following notation: Given $c_1,\ldots,c_k\in C$, we denote  $c_1\cdot\cdot\cdot c_k=\mu_{2,0}(\ldots\mu_{2,0}(\mu_{2,0}(c_1,c_2),c_3)$ $\ldots,c_k)$ with the convention that if $k=0$ we mean $e_C$.

\begin{prop}\label{maptensor}
Let $C_1, C_2, D_1$ and $D_2$ be filtered $dg$-algebras. Let $\varphi, f: C_1\lto C_2$ and $\psi, g:D_1\lto D_2$ be unital filtered $A_\infty$-homomorphisms with $f$ and $g$ naive. Then there are \Ai-homomorphisms $\varphi\otimes_{dg} g , \ f \otimes_{dg} \psi:C_1\otimes_{dg}D_1\lto C_2\otimes_{dg}D_2$ defined as follows:
 $$(\varphi\otimes_{dg} g)_{k,\beta_1} = (-1)^{\lhd} \varphi_{k,\beta_1}(c_1, \ldots, c_k) \otimes g(d_1\cdot\cdot\cdot d_k),$$
 if $\beta_1 \in G_C$ and zero otherwise;
 $$ (f \otimes_{dg} \psi)_{k,\beta_2} = (-1)^{\lhd + (k+1)\vert {\bf c } \vert} f(c_1\cdot\cdot\cdot c_k) \otimes \psi_{k,\beta_2}(d_1,\ldots,d_k)$$
 if $\beta_2 \in G_D$ and zero otherwise.
 
 In either case, the sign is defined in terms of $\lhd:=\sum_{p>q}\vert c_p\vert \vert d_q\vert$.
\end{prop}
\begin{proof}
Since both statements have essentially the same proof, we only derive the \Ai-homomorphism equation for $\varphi \otimes_{dg} g$.

We start by observing that the \Ai-equations of a filtered dg-algebra imply the following identities:
\begin{align}
&(c_1\cdot\cdot\cdot c_j)(c_{j+1}\cdot\cdot\cdot c_k)=(-1)^{(k-j-1)\sum_{l=1}^j\vert c_l\vert}c_1\cdot\cdot\cdot c_k,\nonumber\\
&c_1\cdot\cdot\cdot(c_i \ c_{i+1})\cdot\cdot\cdot c_k=(-1)^{\sum_{l=1}^{i-1}\vert c_l \vert}c_1\cdot\cdot\cdot c_k,\nonumber\\
&c_1\cdot\cdot\cdot c_i\ e_C\ c_{i+1}\cdot\cdot\cdot c_k=(-1)^{\sum_{l=1}^i\vert c_l\vert}c_1 \cdot\cdot\cdot c_k, \label{dgeq} \\
&(-1)^k\mu_{1,\beta}(c_1\cdot\cdot\cdot c_k)+\sum_{1\leq i \leq k }(-1)^{\sum_{l=1}^{i-1}\vert\vert c_l \vert\vert}c_1\cdot\cdot\cdot\mu_{1,\beta}(c_i)\cdot\cdot\cdot c_k=0.\nonumber
\end{align}
We omit the derivation of the \Ai-homomorphism equation with no inputs as it is very similar to the general case. Thus we will prove that the following equation holds for all $k>0$.
\begin{align}
&\sum_{\substack{0\leq j\leq k\\ \beta'+\beta''=\beta}}\mu^\otimes_{2,0}\big((\varphi\otimes g)_{j,\beta'}(c_1\otimes d_1,\ldots,c_j\otimes d_j),(\varphi\otimes g)_{k-j,\beta''}(c_{j+1}\otimes d_{j+1},\ldots,c_k\otimes d_k)\big)+\nonumber\\
&+\sum_{\beta'+\beta''=\beta}\mu^\otimes_{1,\beta''}\big((\varphi\otimes g)_{k,\beta'}(c_1\otimes d_1,\ldots,c_k\otimes d_k)\big)=\nonumber\\
&=\sum_{1\leq i\leq k-1}(-1)^*(\varphi\otimes g)_{k-1,\beta}(c_1\otimes d_1,\ldots,\mu^\otimes_{2,0}(c_i\otimes d_i,c_{i+1}\otimes d_{i+1}),\ldots,c_k\otimes d_k) \nonumber\\
&+\sum_{\substack{1\leq i\leq k\\\beta'+\beta''=\beta}}(-1)^*(\varphi\otimes g)_{k,\beta'}(c_1\otimes d_1,\ldots,\mu^\otimes_{1,\beta''}(c_i\otimes d_i),\ldots,c_k\otimes d_k)\nonumber\\
&+\sum_{\substack{0\leq i\leq k\\\beta'+\beta''=\beta}}(-1)^*(\varphi\otimes g)_{k+1,\beta'}(c_1\otimes d_1,\ldots,c_i\otimes d_i,\mu_{0,\beta''}^\otimes,c_{i+1}\otimes d_{i+1},\ldots,c_k\otimes d_k),\nonumber
\end{align}
We will do this by breaking this equation into a sum of two equalities. The first consists of the terms where both $\beta',\beta'' \in G_C$ and of these, only those involving $\mu_{k,\beta' \times 0}^\otimes$. The second equation has the remaining terms, namely those involving terms of the form $\mu_{k,0 \times \beta''}^\otimes$. The first one equals
\begin{align}
&\sum_{\beta',\beta'' \in G_C}(-1)^{\lhd +(k-j-1)\sum_{l=1}^j\vert d_l\vert}\mu_{2,0}\big(\varphi_{j,\beta'}(c_1,\ldots,c_j),\varphi_{k-j,\beta''}(c_{j+1},\ldots,c_k)\big)\otimes\nonumber\\
& \ \ \ \ g(d_1\cdot\cdot\cdot d_j) g(d_{j+1}\cdot\cdot\cdot d_k)+(-1)^{\lhd}\sum_{\beta',\beta'' \in G_C}\mu_{1,\beta'}(\varphi_{k,\beta''}(c_1,\ldots,c_k))\otimes g(d_1\cdot\cdot\cdot d_k)\nonumber\\
&=\sum_{1\leq i\leq k-1}(-1)^{\zeta+\sum_{l=1}^{i-1}\vert d_l\vert}\varphi_{k-1,\beta}(c_1,\ldots,\mu_{2,0}(c_i,c_{i+1}),\ldots,c_k)\otimes g(d_1\cdot\cdot\cdot(d_i d_{i+1})\cdot\cdot\cdot d_k)\nonumber\\
&+\sum_{\substack{1\leq i\leq k\\\beta',\beta'' \in G_C}}(-1)^{\zeta}(c_1,\ldots,\mu_{1,\beta''}(c_i),\ldots,c_k)\otimes g(d_1\cdot\cdot\cdot d_k)\nonumber\\
&+\sum_{\substack{\beta',\beta''\in G_C \\0\leq i\leq k}}(-1)^{\zeta+\sum_{l=1}^{i-1}\vert d_l\vert} \varphi_{k+1,\beta'}(c_1,\ldots,\mu_{0,\beta''},\ldots,c_k)\otimes g(d_1\cdot\cdot\cdot d_i\ e_{D_1}\ d_{i+1}\cdot\cdot\cdot d_k),\nonumber
\end{align}
where $\zeta=\lhd+\sum_{l=1}^{i-1}\vert\vert c_l\vert\vert $. Using the first three identities in (\ref{dgeq}) and the fact that $g$ is a naive \Ai-homomorphism we conclude that this equation holds, since it is equivalent to the \Ai-homomorphism equation for $\varphi_{k,\beta}$.

As for the second equation, we have
\begin{align}
&\sum_{\beta'+\beta''=\beta}(-1)^{\lhd +k+1+\vert {\bf c}\vert}\varphi_{k,\beta'}(c_1,\ldots,c_k)\otimes \mu_{1,\beta''}g((d_1\cdot\cdot\cdot d_k))\nonumber\\
&=\sum_{1\leq i\leq k}(-1)^{\lhd +\vert {\bf c}\vert +\sum_{l=1}^{i-1}\vert\vert d_l \vert\vert}\varphi_{k,\beta'}(c_1,\ldots,c_k)\otimes g(d_1\cdot\cdot\cdot\mu_{1, \beta''}(d_i)\cdot\cdot\cdot d_k)\nonumber\\
&+\sum_{0\leq i\leq k}(-1)^{\lhd + \sum_{l=1}^{i-1}(\vert\vert c_l\vert\vert +\vert d_l \vert)}\varphi_{k+1,\beta'}(c_1,\ldots,e_{C_1},\ldots,c_k)\otimes g(d_1 \cdot\cdot\cdot d_i\cdot \mu_{0,\beta''}\cdot\cdot\cdot d_k).\nonumber
\end{align}

For $k>0$, the last sum vanishes because $\varphi$ is unital, therefore the equation is equivalent to
\begin{align}
&\sum_{\beta'+\beta''=\beta}(-1)^{\lhd +\vert{\bf c}\vert}\varphi_{k,\beta'}(c_1,\ldots,c_k)\otimes \nonumber\\
&\otimes\left[(-1)^k\mu_{1,\beta''}(d_1\cdot\cdot\cdot d_k)+\sum_{1\leq i\leq k}(-1)^{\sum_{l=1}^{i-1}\vert\vert d_l\vert\vert}d_1\cdot\cdot\cdot\mu_{1,\beta_2}(d_i)\cdot\cdot\cdot d_k\right]=0.\nonumber
\end{align}
Finally this equation holds because the right hand side of the tensor product vanishes by the last identity in (\ref{dgeq}).
\end{proof}

This proposition allows us to give the following definition:
\begin{defn}
 Let $\varphi: C_1\lto C_2$ and $\psi:D_1\lto D_2$ be filtered, unital \Ai-homomorphisms between filtered dg-algebras $C_1, C_2, D_1,D_2$. Define $\varphi\otimes_{dg}\psi:C_1\otimes_{dg}D_1\lto C_2\otimes_{dg}D_2$ by
$$\varphi\otimes_{dg}\psi = (\varphi\otimes_{dg} id_{D_1})\circ (id_{C_2}\otimes_{dg} \psi).$$
\end{defn}

We would like to point out that this construction is not functorial, that is, 
$$(\varphi_1\otimes_{dg} \psi_1)\circ (\varphi_2\otimes_{dg} \psi_2)\neq (\varphi_1\circ \varphi_2)\otimes_{dg} (\psi_1\circ \psi_2).$$ 
We believe that these two \Ai-homomorphism are homotopy equivalent, but we have not tried to prove this. In fact we will not need this statement, instead we will use the following 

\begin{lem}\label{dgqism}
 Let $\varphi: C_1\lto C_2$ and $\psi:D_1\lto D_2$ be filtered quasi-isomorphisms between filtered dg-algebras $C_1, C_2, D_1$ and $D_2$. Then $\varphi\otimes_{dg}\psi:C_1\otimes_{dg}D_1\lto C_2\otimes_{dg}D_2$ is a quasi-isomorphism.
\end{lem}
\begin{proof}
We only need to check that $(\varphi\otimes_{dg}\psi)_{1,0}$ is an isomorphism on cohomology. For this note that $(\varphi\otimes_{dg}\psi)_{1,0}(c\otimes d)=\varphi_{1,0}(c)\otimes \psi_{1,0}(d)$ by definition. Then the K\"unneth theorem implies that 
$$\varphi_{1,0}\otimes\psi_{1,0}:C_1\otimes D_1\lto C_2\otimes D_2$$
induces an isomorphism in cohomology.
\end{proof}

\medskip

We now tackle the problem for general \Ai-algebras, $A_1$, $A_2$, $B_1$ and $B_2$.
In the previous section we constructed a quasi-isomorphism $F_{A_i}:A_i\lto End_{A_i}$; let $F^{-1}_{A_i}$ be a quasi-inverse. In fact, we could use the homological perturbation lemma to give explicit formulas for $F^{-1}_{A_i}$, but this is unnecessary for our purposes.

Given $\varphi: A_1\lto A_2$ and $\psi: B_1\lto B_2$ define $\tilde\varphi:=F_{A_2}\circ\varphi\circ F_{A_1}^{-1}$ and $\tilde\psi:=F_{B_2}\circ\psi\circ F_{B_1}^{-1}$. Finally, define $\varphi\otimes_\infty\psi:A_1\otimes_\infty B_1\lto A_2\otimes_\infty B_2$ by
$$\varphi\otimes_\infty\psi:= \tilde\varphi\otimes_{dg}\tilde\psi: End_{A_1}\otimes_{dg} End_{B_1}\lto End_{A_2}\otimes_{dg} End_{B_2}.$$

\begin{prop}
 Let $\varphi: A_1\lto A_2$ and $\psi: B_1\lto B_2$ be quasi-isomorphism. Then $\varphi\otimes_\infty\psi:A_1\otimes_\infty B_1\lto A_2\otimes_\infty B_2$ is a quasi-isomorphism
\end{prop}
\begin{proof}
We only need to check that $(\varphi\otimes_\infty\psi)_{1,0}$ is an isomorphism on cohomology.
First note that, if $\varphi$ and $\psi$ are quasi-isomorphisms then so are $\tilde\varphi$ and $\tilde\psi$. Thus the statement follows from Lemma \ref{dgqism} applied to $\tilde\varphi$ and $\tilde\psi$.
\end{proof}

Next we will discuss some monoidal properties of the tensor product of $A_\infty$-algebras. We will not show that the category of filtered $A_\infty$-algebras is monoidal, which is probably true but not relevant here.

\begin{prop} Let $A$, $B$ and $C$ be filtered $A_\infty$-algebras. We have the following:
\begin{itemize}
\item[{\bf(a)}] Consider the field $\mathbb{K}$ as a filtered \Ai-algebra with the product as the single non-trivial operation $\m_{2,0}$. Then $A\otimes_\infty \mathbb{K}\simeq A$.
\item[{\bf(b)}] $A\otimes_\infty B\simeq B\otimes_\infty A$.
\item[{\bf(c)}] $A\otimes_\infty(B\otimes_\infty C)\simeq(A\otimes_\infty B)\otimes_\infty C$.
\end{itemize} 
\end{prop}
\begin{proof}
For part (a) note that, given a filtered dg-algebra $V$, we have $V \otimes_{dg} \mathbb{K} \simeq V$. Also, by definition $End_{\mathbb{K}}$ is isomorphic to $\mathbb{K}$, hence
$$A\otimes_\infty \mathbb{K} = End_A \otimes_{dg} End_{\mathbb{K}} \simeq End_A \otimes_{dg} \mathbb{K} \simeq A,$$
by Lemma \ref{dgqism}.

Given filtered dg-algebras $V$ and $W$ define the naive map $$\tau:V\otimes_{dg} W\lto W\otimes_{dg} V$$ by $\tau_{1,0}(v\otimes w)=(-1)^{\vert v\vert \vert w \vert}w\otimes v$. One can easily see this is a \Ai-isomorphism. Applying this to $End_A$ and $End_B$ proves (b).

Next recall that $End_B\otimes_{dg} End_C\simeq End_{End_B\otimes_{dg} End_C}$, by Corollary \ref{EndA}, therefore Lemma \ref{dgqism} implies
\begin{align}
A\otimes_\infty(B\otimes_\infty C)&:= End_A\otimes_{dg} End_{End_B\otimes_{dg} End_C}\nonumber\\
&\simeq End_A\otimes_{dg}(End_B\otimes_{dg} End_C).\nonumber
\end{align}
Similarly
$$(A\otimes_\infty B)\otimes_\infty C\simeq(End_A\otimes_{dg} End_B)\otimes_{dg} End_C.$$
Part (c) now follows from the following easily verifiable fact: Given filtered dg-algebras $V$, $W$ and $U$, the naive map
$$id_{1,0}:V\otimes_{dg}(W\otimes_{dg} U)\lto (V\otimes_{dg}W)\otimes_{dg} U$$ is an $A_\infty$-isomorphism.
\end{proof}

\section{Another model for $A\otimes_\infty B$}

Given two \Ai-algebras $A$ and $B$, we will use the homological perturbation lemma to construct an \Ai-algebra structure on the vector space $A\otimes B$ quasi-isomorphic to $A\otimes_\infty B$.

We start by introducing some terminology and constructions relating to trees.
Let $G_n$ denote the set of all trees with $n$ leaves such that $val(v)\geq 2$ for all $v\in V(T)$. Given a vertex $v$ we will denote by $v(i)$ the $i^{th}$ incident edge at $v$. 

For each pair $x,y \in V(T) \cup E(T)$ we say $x$ is under $y$ if the unique path $P_x$ in $T$ from $x$ to the root of $T$ contains $y$. Note that this defines a partial ordering. Suppose $x$ and $y$ are not comparable and denote by $w$ the first (furthest from the root) vertex that $P_x$ and $P_y$ have in common. Then $P_x$ and $P_y$ contain edges $w(i)$ and $w(j)$ respectively. We say that $x$ is to the left of $y$ if $i<j$.

Given $x \in V(T) \cup E(T)$ we denote by $p_x$ ( respectively $q_x$) the left  (respectively right) most leaf of $T$ under $x$. Additionally, given $e \in E(T)$ we define integers $J_e$, $I_e$ and $\tilde{I_e}$ to be zero when $e$ is a leaf and when $e$ is not a leaf, $J_e$ is the number of leaves under $e$, $\tilde{I_e}$ is the number of internal edges under $e$ and $I_e = \tilde{I_e} +1$.

Next, let $G^{bin}_n$ denote the subset of $G_n$ consisting of trees $T$ with $val(v)= 2$ for all $v\in V(T)$. Elements of $G^{bin}_n$ are called binary trees. There is a partial order on the set $G^{bin}_n$, which was introduced in \cite{MarShn}. It is generated by the relation on the following figure:
\begin{figure}[h]
\begin{center}
 \setlength{\unitlength}{3pt}
\begin{picture}(60,20)(0,0)

\linethickness{0.3mm}
\put(10,10){\line(0,1){10}}
\linethickness{0.3mm}
\put(0,0){\line(1,1){10}}
\linethickness{0.3mm}
\put(20,0){\line(-1,1){10}}
\linethickness{0.3mm}
\put(10,0){\line(-1,1){5}}

\linethickness{0.3mm}
\put(50,10){\line(0,1){10}}
\linethickness{0.3mm}
\put(40,0){\line(1,1){10}}
\linethickness{0.3mm}
\put(60,0){\line(-1,1){10}}
\linethickness{0.3mm}
\put(50,0){\line(1,1){5}}

\put(30,10){\makebox(0,0)[cc]{$\leq$}}

\end{picture}
\end{center}
\caption{Relation that generates the partial order on binary trees}
\label{figordertree}
\end{figure}
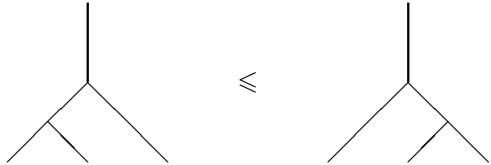

We say $T_1\leq T_2$ if there is a sequence of trees starting at $T_1$ and ending at $T_2$ such that consecutive trees are isomorphic in the complement of some neighborhood where they differ as pictured in Figure \ref{figordertree}. There is an absolute minimum and maximum for this order. The minimum (maximum), which we denote by $\underline{M}_n$ ($\overline{M}_n$), is the binary tree with strictly right (left) leaning internal edges.

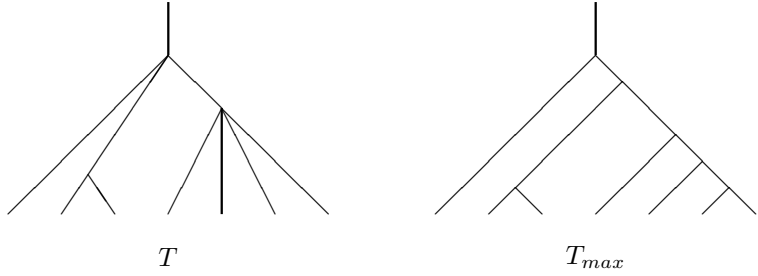
\begin{figure}[h]
\begin{center}
 \setlength{\unitlength}{2pt}
\begin{picture}(140,50)(0,0)

\linethickness{0.3mm}
\put(30,40){\line(0,1){10}}
\linethickness{0.3mm}
\put(0,10){\line(1,1){30}}
\linethickness{0.3mm}
\put(10,10){\line(2,3){20}}
\linethickness{0.3mm}
\put(20,10){\line(-2,3){5}}
\linethickness{0.3mm}
\put(30,10){\line(1,2){10}}
\linethickness{0.3mm}
\put(40,10){\line(0,1){20}}
\linethickness{0.3mm}
\put(50,10){\line(-1,2){10}}
\linethickness{0.3mm}
\put(60,10){\line(-1,1){30}}

\linethickness{0.3mm}
\put(110,40){\line(0,1){10}}
\linethickness{0.3mm}
\put(80,10){\line(1,1){30}}
\linethickness{0.3mm}
\put(90,10){\line(1,1){25}}
\linethickness{0.3mm}
\put(100,10){\line(-1,1){5}}
\linethickness{0.3mm}
\put(110,10){\line(1,1){15}}
\linethickness{0.3mm}
\put(120,10){\line(1,1){10}}
\linethickness{0.3mm}
\put(130,10){\line(1,1){5}}
\linethickness{0.3mm}
\put(140,10){\line(-1,1){30}}

\put(30,2){\makebox(0,0)[cc]{$T$}}
\put(110,2){\makebox(0,0)[cc]{$T_{max}$}}

\end{picture}
\end{center}
\caption{Example of maximal tree resolving $T$}
\label{figmaxtree}
\end{figure}

Given $T\in G_n$ we define $T_{max}\in G^{bin}_n$ as the maximal (with respect to the order just described) binary tree that resolves $T$. That is, $T$ can be obtained by collapsing several edges in $T_{max}$. See Figure \ref{figmaxtree} for an example.

Similarly, we define $T_{min}$ as the minimal binary tree that resolves $T$.
The tree in $G_n$ with exactly one internal vertex, called the \emph{$n$-corolla}, is denote by $c_n$.

We introduce one more subset of $G_n$. Let $L_n$ denote the subset of $G_n$ obtained by grafting corollas in all possible ways while avoiding the last leaf. More formally, $T\in L_n$ if either $T=c_n$ or there exist integers $1\leq i_1 < \ldots < i_l < m$ and trees $T_j \in L_{n_j}$ such that
$$T=c_m\circ_{i_1,\ldots,i_l}(T_1,\ldots ,T_l ),$$ 
where $\circ_{i_1,\ldots,i_l}$ tells us to graft the root of $T_j$ to the leaf $i_j$ of $c_m$.

\begin{figure}[h]
\begin{center}
 \setlength{\unitlength}{1.8pt}
\begin{picture}(130,30)(0,0)

\put(3,10){\makebox(0,0)[cc]{$\alpha\Bigg($}}

\linethickness{0.3mm}
\put(25,20){\line(0,1){10}}
\linethickness{0.3mm}
\put(10,0){\line(3,4){15}}
\linethickness{0.3mm}
\put(20,0){\line(0,1){13,3}}
\linethickness{0.3mm}
\put(30,0){\line(-3,4){10}}
\linethickness{0.3mm}
\put(40,0){\line(-3,4){15}}
\linethickness{0.3mm}

\put(45,10){\makebox(0,0)[cc]{$\Bigg)$}}
\put(55,10){\makebox(0,0)[cc]{$=$}}

\linethickness{0.3mm}
\put(75,10){\line(0,1){20}}
\linethickness{0.3mm}
\put(60,0){\line(3,4){15}}
\linethickness{0.3mm}
\put(70,0){\line(1,2){5}}
\linethickness{0.3mm}
\put(80,0){\line(-1,2){5}}
\linethickness{0.3mm}
\put(90,0){\line(-3,4){15}}
\linethickness{0.3mm}

\put(95,10){\makebox(0,0)[cc]{$+$}}

\linethickness{0.3mm}
\put(115,20){\line(0,1){10}}
\linethickness{0.3mm}
\put(100,0){\line(3,4){15}}
\linethickness{0.3mm}
\put(110,0){\line(3,4){10}}
\linethickness{0.3mm}
\put(120,0){\line(0,1){13,3}}
\linethickness{0.3mm}
\put(130,0){\line(-3,4){15}}
\linethickness{0.3mm}

\put(130,10){\makebox(0,0)[cc]{$.$}}
\end{picture}
\end{center}
\caption{Example of $\alpha$} 
\label{figexalpha}
\end{figure}
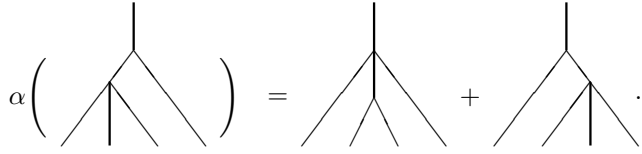

Now suppose that $V\in L_n$ has $k$ internal edges. There is an obvious correspondence between internal edges of $V$ and right leaning (internal) edges of $V_{max}$. Let $R(B)$ denote the set of right  leaning internal edges of a binary tree $B$ and let $\vert R(B)\vert$ denote its order. We define
$$\alpha(V):=\sum_{\substack{S\in G^{bin}_n\\S\geq V_{max}\\ \vert R(S)\vert=k}}S/R(S),$$
where $S/R(S)$ is the tree obtained by collapsing all the edges of $S$ in $R(S)$. Note that $\alpha(V)$ is non-zero as it always includes $\alpha_0(T):= T_{max}/R(T_{max})$. For an example see Figure \ref{figexalpha}.

Finally, given $U \in G_n$ and $\beta \in G$, let $U_\beta$ denote the sum of all elements of $\Gamma_n(\beta)$ whose underlying tree is $U$. We use $U_\beta$ as a flow chart to define a map
$$U_\beta^A: A^{\otimes k}\lto A$$
by assigning to each $v\in V(U_\beta)$ the map $\m^A_{val(v),\nu(v)}$. Similarly we define maps $U_\beta^B$.

\begin{thm}\label{filteredtensor}
Let $A$ and $B$ be unital filtered $A_\infty$-algebras. Then  the tensor product $A\otimes_\infty B $ is quasi-isomorphic to the $A_\infty$-algebra $( A\otimes B, \m^{\otimes})$ with operations given by $\m_{k,\beta}^{\otimes}=\sum_{\beta_1+\beta_2=\beta}\m_{k,\beta_1 \times \beta_2}^{\otimes}$, where
\begin{align}
\m^{\otimes}_{0,\beta_1\times\beta_2}&=\left\{\begin{array}{ll}
\m^{A}_{0,\beta_1}\otimes e_B&,\beta_2=0\\
e_A\otimes \m_{0,\beta_2}^B&,\beta_1=0\\
0&, \beta_1, \beta_2 \neq 0
\end{array},\right.\nonumber\end{align}
\begin{align}
\m^{\otimes}_{1,\beta_1\times\beta_2}&=\left\{\begin{array}{ll}
\m_{1,\beta_1}\otimes id&,\beta_2=0\\
id\otimes \m_{1,\beta_2}&,\beta_1=0\\
0&, \beta_1, \beta_2 \neq 0
\end{array},\right. \label{tensor} \\
\m^{\otimes}_{n,\beta_1\times\beta_2}(a_1 \otimes b_1, \ldots , a_n \otimes b_n)&=\sum_{T\in L_n} (-1)^{\epsilon} \ T^A_{\beta_1} (a_1, \ldots a_n) \otimes \alpha(T)^B_{\beta_2}(b_1,\ldots,b_n). \nonumber
\end{align}
The sign exponent is defined as $\epsilon= \theta(T) + \lhd + \vert {\bf a}\vert \vert E_{int}(T) \vert + \gamma_{\bf a} + \gamma_{\bf b}$, with
$$\gamma_{\bf a}= \sum_{v \in V(T)} \sum_{i< p_v} \vert a_i \vert, \ \ 
\gamma_{\bf b}= \sum_{v \in V(\alpha(T))} \sum_{i< p_v} \vert b_i \vert \ \textrm{and}$$
$$\theta(T)= \sum_{v \in V(T)} \big( \sum_{i} ((val(v)-i-1)I_{v(i)} + \tilde{I}_{v(i)} +J_{v(i)}) + \sum_{i<j} I_{v(i)}(\tilde{I}_{v(j)} + J_{v(j)})\big).$$

\end{thm}

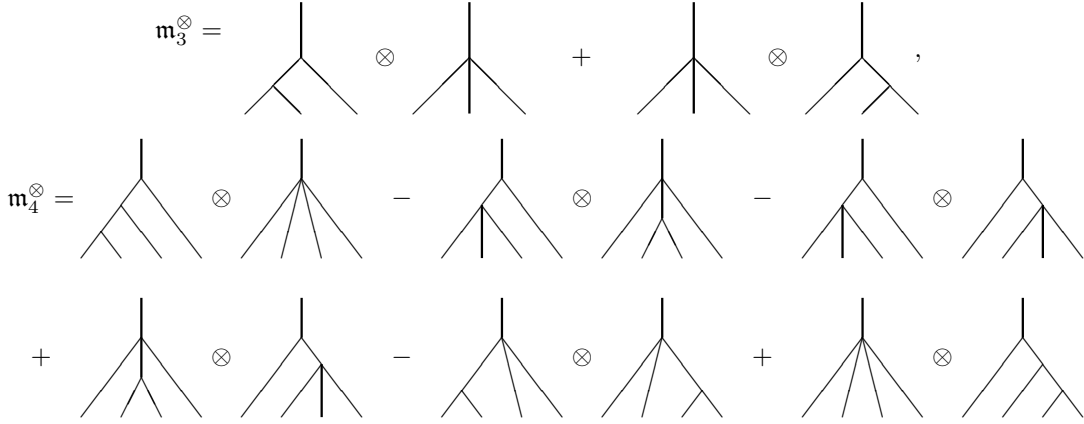
\begin{figure}[h]
\begin{center}
 \setlength{\unitlength}{2.1pt}
\begin{picture}(130,20)(0,0)
\put(0,15){\makebox(0,0)[cc]{$\m_3^{\otimes}=$}}

\linethickness{0.3mm}
\put(20,10){\line(0,1){10}}
\linethickness{0.3mm}
\put(10,0){\line(1,1){10}}
\linethickness{0.3mm}
\put(30,0){\line(-1,1){10}}
\linethickness{0.3mm}
\put(20,0){\line(-1,1){5}}

\put(35,10){\makebox(0,0)[cc]{$\otimes$}}

\linethickness{0.3mm}
\put(50,0){\line(0,1){20}}
\linethickness{0.3mm}
\put(40,0){\line(1,1){10}}
\linethickness{0.3mm}
\put(60,0){\line(-1,1){10}}

\put(70,10){\makebox(0,0)[cc]{$+$}}

\linethickness{0.3mm}
\put(90,0){\line(0,1){20}}
\linethickness{0.3mm}
\put(80,0){\line(1,1){10}}
\linethickness{0.3mm}
\put(100,0){\line(-1,1){10}}

\put(105,10){\makebox(0,0)[cc]{$\otimes$}}

\linethickness{0.3mm}
\put(120,10){\line(0,1){10}}
\linethickness{0.3mm}
\put(110,0){\line(1,1){10}}
\linethickness{0.3mm}
\put(130,0){\line(-1,1){10}}
\linethickness{0.3mm}
\put(120,0){\line(1,1){5}}

\put(130,10){\makebox(0,0)[cc]{$,$}}

\end{picture}
\vspace{0.3cm}
%%%%%

 \setlength{\unitlength}{1.5pt}
\begin{picture}(270,70)(0,0)
\put(0,55){\makebox(0,0)[cc]{$\m_4^{\otimes}=$}}

\linethickness{0.3mm}
\put(25,60){\line(0,1){10}}
\linethickness{0.3mm}
\put(10,40){\line(3,4){15}}
\linethickness{0.3mm}
\put(40,40){\line(-3,4){15}}
\linethickness{0.3mm}
\put(20,40){\line(-3,4){5}}
\linethickness{0.3mm}
\put(30,40){\line(-3,4){10}}

\put(45,55){\makebox(0,0)[cc]{$\otimes$}}

\linethickness{0.3mm}
\put(65,60){\line(0,1){10}}
\linethickness{0.3mm}
\put(50,40){\line(3,4){15}}
\linethickness{0.3mm}
\put(80,40){\line(-3,4){15}}
\linethickness{0.3mm}
\put(60,40){\line(1,4){5}}
\linethickness{0.3mm}
\put(70,40){\line(-1,4){5}}

\put(90,55){\makebox(0,0)[cc]{$-$}}

\linethickness{0.3mm}
\put(115,60){\line(0,1){10}}
\linethickness{0.3mm}
\put(100,40){\line(3,4){15}}
\linethickness{0.3mm}
\put(130,40){\line(-3,4){15}}
\linethickness{0.3mm}
\put(110,40){\line(0,1){13,33}}
\linethickness{0.3mm}
\put(120,40){\line(-3,4){10}}

\put(135,55){\makebox(0,0)[cc]{$\otimes$}}

\linethickness{0.3mm}
\put(155,50){\line(0,1){20}}
\linethickness{0.3mm}
\put(140,40){\line(3,4){15}}
\linethickness{0.3mm}
\put(170,40){\line(-3,4){15}}
\linethickness{0.3mm}
\put(150,40){\line(1,2){5}}
\linethickness{0.3mm}
\put(160,40){\line(-1,2){5}}

\put(180,55){\makebox(0,0)[cc]{$-$}}

\linethickness{0.3mm}
\put(205,60){\line(0,1){10}}
\linethickness{0.3mm}
\put(190,40){\line(3,4){15}}
\linethickness{0.3mm}
\put(220,40){\line(-3,4){15}}
\linethickness{0.3mm}
\put(200,40){\line(0,1){13,33}}
\linethickness{0.3mm}
\put(210,40){\line(-3,4){10}}

\put(225,55){\makebox(0,0)[cc]{$\otimes$}}

\linethickness{0.3mm}
\put(245,60){\line(0,1){10}}
\linethickness{0.3mm}
\put(230,40){\line(3,4){15}}
\linethickness{0.3mm}
\put(260,40){\line(-3,4){15}}
\linethickness{0.3mm}
\put(240,40){\line(3,4){10}}
\linethickness{0.3mm}
\put(250,40){\line(0,1){13,33}}

%\put(270,55){\makebox(0,0)[cc]{$+$}}

\put(0,15){\makebox(0,0)[cc]{$+$}}

\linethickness{0.3mm}
\put(25,10){\line(0,1){20}}
\linethickness{0.3mm}
\put(10,0){\line(3,4){15}}
\linethickness{0.3mm}
\put(40,0){\line(-3,4){15}}
\linethickness{0.3mm}
\put(20,0){\line(1,2){5}}
\linethickness{0.3mm}
\put(30,0){\line(-1,2){5}}

\put(45,15){\makebox(0,0)[cc]{$\otimes$}}

\linethickness{0.3mm}
\put(65,20){\line(0,1){10}}
\linethickness{0.3mm}
\put(50,0){\line(3,4){15}}
\linethickness{0.3mm}
\put(80,0){\line(-3,4){15}}
\linethickness{0.3mm}
\put(60,0){\line(3,4){10}}
\linethickness{0.3mm}
\put(70,0){\line(0,1){13,33}}

\put(90,15){\makebox(0,0)[cc]{$-$}}

\linethickness{0.3mm}
\put(115,20){\line(0,1){10}}
\linethickness{0.3mm}
\put(100,0){\line(3,4){15}}
\linethickness{0.3mm}
\put(130,0){\line(-3,4){15}}
\linethickness{0.3mm}
\put(110,0){\line(-3,4){5}}
\linethickness{0.3mm}
\put(120,0){\line(-1,4){5}}

\put(135,15){\makebox(0,0)[cc]{$\otimes$}}

\linethickness{0.3mm}
\put(155,20){\line(0,1){10}}
\linethickness{0.3mm}
\put(140,0){\line(3,4){15}}
\linethickness{0.3mm}
\put(170,0){\line(-3,4){15}}
\linethickness{0.3mm}
\put(150,0){\line(1,4){5}}
\linethickness{0.3mm}
\put(160,0){\line(3,4){5}}

\put(180,15){\makebox(0,0)[cc]{$+$}}

\linethickness{0.3mm}
\put(205,20){\line(0,1){10}}
\linethickness{0.3mm}
\put(190,0){\line(3,4){15}}
\linethickness{0.3mm}
\put(220,0){\line(-3,4){15}}
\linethickness{0.3mm}
\put(200,0){\line(1,4){5}}
\linethickness{0.3mm}
\put(210,0){\line(-1,4){5}}

\put(225,15){\makebox(0,0)[cc]{$\otimes$}}

\linethickness{0.3mm}
\put(245,20){\line(0,1){10}}
\linethickness{0.3mm}
\put(230,0){\line(3,4){15}}
\linethickness{0.3mm}
\put(260,0){\line(-3,4){15}}
\linethickness{0.3mm}
\put(240,0){\line(3,4){10}}
\linethickness{0.3mm}
\put(250,0){\line(3,4){5}}

\put(270,15){\makebox(0,0)[cc]{$.$}}

\end{picture}
\end{center}
\caption{Formulas for $\m_3^{\otimes}$ and $\m_4^{\otimes}$.}
\label{figm3m4}
\end{figure}

To illustrate the theorem, explicit formulas for $\m_3^\otimes$ and $\m_4^\otimes$ appear in Figure \ref{figm3m4} with signs determined by $\theta(T)$.

%\vspace{0.3cm}

Assume $A$ and $B$ are classical \Ai-algebras, then (\ref{tensor}) agrees exactly with the formula for $\m_n^{\otimes_{SU}}$ provided in \cite[Theorem 3]{MarShn}, up to sign. The difference in the sign comes from the fact that Markl and Shnider use a different sign rule in the definition of \Ai-algebra. In (\ref{Ainf}), they take $*=i(j+1)+jn+j(\vert a_1\vert+\ldots +\vert a_{i-1} \vert)$. This definition is equivalent to ours via the change 
$$\m'_k (a_1,\ldots,a_n)\mapsto (-1)^{\sum (k-i)\vert a_i \vert}\m_k(a_1,\ldots,a_n).$$
Taking this sign change in consideration one can show that (\ref{tensor}) agrees with the formula in \cite{MarShn}. Therefore we obtain the following

\begin{cor}Let $A$ and $B$ be classical \Ai-algebras, then the tensor product $A\otimes_\infty B$ is quasi-isomorphic to the tensor product $A\otimes_{SU}B$ defined by Markl and Shnider in \cite{MarShn}.
\end{cor}

Before proving the theorem, we will give two alternative descriptions of $\m_{n,\beta_1 \times \beta_2}^{\otimes}$ for $n\geq 3$.
Let $R_n$ denote the subset of $G_n$ obtained by grafting corollas in all possible ways while avoiding the first leaf. Given $T\in R_n$ with $k$ internal edges we define 
$$\bar \alpha(T)=\sum_{\substack{W\in G_n^{bin}\\W\leq T_{min}\\ \vert L(W)\vert =k}}W/L(W),$$ where $W/L(W)$ is obtained by collapsing all the edges in $L(W)$, the set of left leaning edges on $W$. We have the following:
\begin{prop}\label{alttensor}
There are the following alternative descriptions when $n\geq 3$
$$\m^\otimes_{n,\beta_1\times\beta_2}=\sum_{T\in R_n} \pm \bar \alpha(T)^A_{\beta_1}\otimes T^B_{\beta_2},$$ and
$$\m^\otimes_{n,\beta_1\times\beta_2}=\sum_{\substack{\vert E_{int}(U)\vert +\vert E_{int}(W)\vert =n-2\\U_{max}\leq W_{min}}}\pm U^A_{\beta_1}\otimes W^B_{\beta_2}.$$
\end{prop}
\begin{proof}
Easy combinatorial check, see \cite{MarShn}.
\end{proof}

The proof of Theorem \ref{filteredtensor} will occupy the remainder of this section.

\begin{proof}[Proof of Theorem \ref{filteredtensor}] This is a direct application of the homological perturbation lemma. Consider the chain complex $(A\otimes B, \delta=\m^A_{1,0}\otimes id + id\otimes \m^B_{1,0})$ and define maps
\begin{align}
 i=F^{A}_{1,0}\otimes F^{A}_{1,0}&: A\otimes B \lto End_A \otimes End_B, \nonumber \\
 P=P_A\otimes P_B &: End_A \otimes End_B \lto A\otimes B, \\
 H= \Pi_A \otimes H_B + H_A \otimes id &: End_A \otimes End_B \lto End_A \otimes End_B,\nonumber   
\end{align}
where $\Pi_A = F^A_{1,0} \circ P_A$ and $P_A$, $H_A$ are as defined in Lemma \ref{homotopytensor}. 

\begin{lem}
 The maps $i$, $P$ and $H$ determine homotopy data that satisfies the side conditions.
\end{lem}
\begin{proof}
The maps $i$ and $P$ are chain maps because they are the tensor product of chain maps. Next, by direct computation we have
\begin{align}
\mu_{1,0}^{\otimes}H+H\mu_{1,0}^{\otimes}&=\Pi_A\otimes\mu_{1,0}^BH_B-H_A\otimes \mu_{1,0}^B+\mu_{1,0}^A\Pi_A\otimes H_B+\mu_{1,0}^A H_A\otimes id+\nonumber\\
+&\Pi_A\otimes H_B\mu_{1,0}^B+H_A\otimes\mu_{1,0}^B-\Pi_A\mu_{1,0}^A\otimes H_B+H_A\mu_{1,0}^A\otimes id\nonumber\\
=&\Pi_A\otimes(\Pi_B-id)+(\mu_{1,0}^A\Pi_A-\Pi_A\mu_{1,0}^A)\otimes H_B+(\Pi_A-id)\otimes id\nonumber\\
=&\Pi_A \otimes \Pi_B - id = i\circ P - id,\nonumber
\end{align}
using Lemma \ref{homotopytensor} for $A$ and $B$ and the fact that $\Pi_A$ is a chain map.

The side conditions follow from direct computations that we omit. Finally we observe that $i(e_A \otimes e_B)=id$.
\end{proof}

Now Theorem \ref{hpl}, provides an \Ai-algebra structure on $A\otimes B$ with $\eta_{1,0}=\delta$. By Proposition \ref{sideconditions} this \Ai-algebra is unital and there is a quasi-isomorphism $\varphi: (A\otimes B, \nu) \lto A \otimes_{\infty}B$. Moreover we have explicit formulas
$$\eta_{n,\beta}=\sum_{T\in\Gamma_n(\beta)} \eta_{T}.$$
In this situation, since the \Ai-algebras are bi-gapped we can further decompose the expression for $\eta_{n,\beta}$ as
\begin{align}\label{formulahpl}
 \eta_{n,\beta}=\sum_{\beta_1+\beta_2=\beta}\eta_{n,\beta \times \beta_2}=\sum_{T\in\Gamma_n(\beta_1\times \beta_2)} \eta_{T}.
\end{align}

All that is left to prove is that (\ref{formulahpl}) gives the formulas in (\ref{tensor}). This will be done in several lemmas in which we will determine and describe the trees $T \in \Gamma_n(\beta_1\times \beta_2)$ that give non-trivial contributions. The first lemma follows from easy computations, which we omit.

\begin{lem}\label{vanishing}
 Let $(C, \m)$ be a $G$-gapped filtered \Ai-algebra and consider $\beta \in G$ and $f, g \in End_C$. We have the following identities:
\begin{itemize}
\item[{\bf(a)}] $H_C(\mu_{0,\beta})=0$.
\item[{\bf(b)}] $H_C\circ\mu_{2,0}(f,H_C (g))=0$, $P_C\circ\mu_{2,0}(f,H_C(g))=0$.
\item[{\bf(c)}] $H_C\circ \mu_{1,\beta}\circ H_C=0$, $P_C\circ \mu_{1,\beta}\circ H_C=0$, if $\beta\neq 0$.
\item[{\bf(d)}] $[H_C\circ \mu_{1,\beta}\circ F^C_{1,0}(\xi)]_*(u,c_1,\ldots,c_r)=\m_{r+2,\beta}(\xi,u,c_1,\ldots,c_r)$.
\item[{\bf(e)}] $P_C\circ \mu_{1,\beta}\circ F^C_{1,0}(\xi)=\m_{1,\beta}(\xi)$.
\end{itemize}
\end{lem}

\begin{lem}\label{m0}
 $H(\mu^{\otimes}_{0,\beta_1\times 0})=H(\mu^{\otimes}_{0,0\times \beta_2})=0$. Therefore $\eta_T=0$ whenever $T$ has a vertex of valency zero.
\end{lem}
\begin{proof}
From the definitions we have
\begin{align}
 H(\mu^{\otimes}_{0,\beta_1\times 0})&= \Pi_A(\mu_{0,\beta_1})\otimes H_B(id_{End_B}) + H_A(\mu_{0,\beta_1})\otimes id_{End_B}=0\nonumber
\end{align}
by Lemma \ref{vanishing}(a) and the easy fact $H_B(id_{End_B})=0$.

The same argument shows $H(\mu^{\otimes}_{0,0\times \beta_2})=0$.

For the second statement, recall that when defining $\eta_T$, a vertex of valency zero is assigned $\mu_{0,\beta}$ and a new vertex inserted into the outgoing edge and assigned $H$. Since we just showed that composition is zero, it follows that $\eta_T=0$.
\end{proof}

This lemma is enough to describe $\eta_{0,\beta}$. A tree $T\in \Gamma_0(\beta_1\times \beta_2)$ such that $\eta_T \neq 0$ cannot have internal edges by Lemma \ref{m0}. Therefore $T$ has to be a tree with one single internal vertex of weight $\beta_1 \times 0$ or $0\times \beta_2$, in these cases we have
$$\eta_{0,\beta_1\times 0}=P(\mu_{0,\beta_1\times 0})=P_A(\mu_{0,\beta_1})\otimes P_B(Id_B)=\m_{0,\beta_1}\otimes e_B.$$ 
and similarly $\eta_{0,0\times\beta_2}=e_A\otimes \m_{0,\beta_2}.$

Let us now consider the general case, $\eta_{n,\beta}=\sum_{T\in\Gamma_n(\beta)} \eta_T$. Since $A\otimes_{\infty}B$ is a filtered dg-algebra, the vertices of $T$ have valency $\leq2$. By Lemma \ref{m0} $T$ cannot have vertices of valency zero. If a vertex has valency one, we say the incoming edge is vertical, if it has valency two, we say the first incoming edge is right leaning and the second is left leaning.

To better describe the maps $\eta_T$, we will further decompose $\eta_T$ into the following sum
$$\eta_{n,\beta}=\sum_{T\in\Gamma_n(\beta)}\eta_T=\sum_{\substack{T\in\Gamma_n(\beta)\\T_H}}\eta_{T_H}.$$
Here $\eta_{T_H}$ is defined in the same way as $\eta_T$ except we assign either $\Pi_A\otimes H_B$ or $H_A\otimes id$ to each internal edge of $T$.  We then sum over the $2^{\vert E_{int}(T)\vert}$ ways of doing this.

\begin{lem}\label{cut}
Given $T\in \Gamma_n(\beta_1\times\beta_2)$ and $T_H$ as above, cut $T$ along the edges with $\Pi_A\otimes H_B$ assigned and let $W_1,\ldots,W_{k+1}$ denote the trees obtained. If $\eta_{T_H}\neq 0$, then each of the $W_i$ is equal to one of the following:
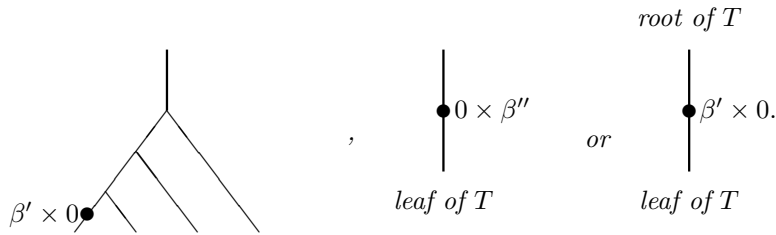
\begin{center}
 \setlength{\unitlength}{2.3pt}
\begin{picture}(120,40)(0,0)

\linethickness{0.3mm}
\put(25,20){\line(0,1){10}}
\linethickness{0.3mm}
\put(10,0){\line(3,4){15}}
\linethickness{0.3mm}
\put(40,0){\line(-3,4){15}}
\linethickness{0.3mm}
\put(20,0){\line(-3,4){5}}
\linethickness{0.3mm}
\put(30,0){\line(-3,4){10}}

\put(12,3){\circle*{2}}

\linethickness{0.3mm}
\put(70,10){\line(0,1){20}}

\put(70,20){\circle*{2}}

\put(70,5){\makebox(0,0)[cc]{leaf of $T$}}
\put(110,5){\makebox(0,0)[cc]{leaf of $T$}}
\put(110,35){\makebox(0,0)[cc]{root of $T$}}

\put(55,15){\makebox(0,0)[cc]{,}}
\put(95,15){\makebox(0,0)[cc]{or}}

\put(5,3){\makebox(0,0)[cc]{$\beta'\times 0$}}
\put(78,20){\makebox(0,0)[cc]{$0\times\beta''$}}

\linethickness{0.3mm}
\put(110,10){\line(0,1){20}}

\put(110,20){\circle*{2}}

\put(118,20){\makebox(0,0)[cc]{$\beta'\times 0.$}}

\end{picture}
\end{center}
In the first case there are $l\geq2$ leaves; if $\beta'=0$, we remove the vertex of valency one. In the second and third cases we have $\beta',\beta'' \neq 0$.
\end{lem}
\begin{proof}
Let $W$ be one such tree. The root of $W$ either coincides with the root of $T$, in which case it is assigned with $P$, or it is assigned with $\Pi_A\otimes H_B$. Let $\eta_W$ denote the map induced by $W$. We are interested in computing $\eta_W(\xi_1,\ldots,\xi_k)$ for $\xi_i=f_i\otimes g_i$, where each $f_i\otimes g_i=\Pi_A (\tilde f_i)\otimes H_B(\tilde g_i)$ or $f_i\otimes g_i=i(a_i\otimes b_i)$ since each leaf of $W$ is either a leaf of $T_H$ or an internal edge of $T_H$ labeled by $\Pi_A\otimes H_B$. By construction, all the internal edges of $W$ are labeled with $H_A\otimes id$.

\underline{Claim}: All the internal edges of $W$ are right leaning.
\begin{itemize}
\item Suppose there is a left leaning internal edge. If it is not adjacent to the top vertex of $W$, it leads to a composition of the form
\begin{align}
 (H_A\otimes id)\circ \mu^\otimes_{2,0}(\rho \otimes \tau, H_A(\varphi)\otimes\psi)=H_A(\mu_{2,0}(\rho,H_A(\varphi)))\otimes \mu_{2,0}(\tau,\psi) =0,\nonumber
\end{align}
which vanishes by Lemma \ref{vanishing}(b). If the edge is adjacent to the top vertex of $W$, we have one of the compositions: $(\Pi_A\otimes H_B)\circ \mu^\otimes_{2,0}(\rho \otimes \tau, H_A(\varphi)\otimes \psi)$ or $(P_A\otimes P_B)\circ \mu^\otimes_{2,0}(\rho \otimes \tau, H_A(\varphi)\otimes \psi)$, depending on whether or not the root of $W$ coincides with the root of $T$. In either case we get zero by Lemma \ref{vanishing}(b).
\item If there is a vertical interior edge, the endpoint of that edge has weight $\beta' \times 0$ or $0 \times \beta''$. In the first case we have one of the compositions $H_A\circ\mu_{1,\beta'}\circ H_A$ or $P_A\circ \mu_{1,\beta'}\circ H_A$, which are both zero by Lemma \ref{vanishing}(c). The second case results in compositions of the form $H_B \circ H_B$ or $P_B \circ H_B$, which are both zero by the side conditions. 
\end{itemize}
This claim immediately implies that $W$ can have at most one vertex of valency one, which must be adjacent to the first leaf of $W$. Moreover, when $W$ has $k>1$ leaves, $T$ is the minimal binary tree, once we forget the vertex of valency one. When $k=1$ then $W$ has a single internal vertex. 

Next we observe that $W=T$ whenever the vertex of valency one has weight $\beta'\times 0$ and $k=1$. If the leaf of $W$ is not a leaf of $T$, the input has the form $\Pi_A (\tilde f_i)\otimes H_B(\tilde g_i)$, which leads to one of the compositions $H_B \circ H_B$ or $P_B \circ H_B$ that vanish as a consequence of the side conditions. Similarly, if the root of $W$ does not coincide with the root of $T$, we get $(\Pi_A\otimes H_B)(\mu_{1,\beta'}\otimes id) \circ i$, which again vanishes by the side conditions. 

Now consider the case when the vertex of valency one has weight $0 \times \beta''$. If the leaf adjacent to this vertex is not a leaf of $T$, then either $(\Pi_A\otimes H_B)(id\otimes \mu_{1,\beta''})(\Pi_A\otimes H_B)$ or $(P_A \otimes P_B)(id\otimes \mu_{1,\beta''})(\Pi_A\otimes H_B)$. Again, the side conditions imply that both vanish. 
%Finally, we prove that this vertex must be adjacent to the root of $W$. If this was not the case, the outgoing edge would be an internal edge of $W$, which is labeled by $H_A\otimes id$, nut this leads to $(H_A\otimes id)(id \otimes \mu_{1,\beta''})\circ i=0$. 
\end{proof}

We are now ready to describe $\eta_{1,\beta_1\times\beta_2}$. Consider $T\in \Gamma_1(\beta_1\times \beta_2)$ and let $T_H$  be as above. Then $T_H$ must be obtained by gluing trees of the second and third types in Lemma \ref{cut}. But it must equal one of these since the leaves of these trees must also be leaves of $T$. Therefore
\begin{align}
\eta_{1,\beta_1\times 0}(a\otimes b)=(p_A\otimes p_B)\circ (\mu_{1,\beta_1}\otimes id)\circ(i(a)\otimes i(b))
=\m_{1,\beta_1}(a)\otimes b\nonumber
\end{align}
by Lemma \ref{vanishing}(e). Similarly $\eta_{1,0\times \beta_2}(a\otimes b)=(-1)^{\vert a\vert}a\otimes \m_{1,\beta_2}(b)$ and $\eta_{1,\beta_1 \times \beta_2}=0$ when $\beta_1, \beta_2 \neq 0$. 

To describe $\eta_{n,\beta}=\sum_{T\in\Gamma_n(\beta)} \eta_T$, for $n\geq 2$ we need two more lemmas.

\begin{lem}\label{lastleaf}
 Let $T \in \Gamma_n(\beta_1 \times \beta_2)$ for $n\geq 2$ and consider the decomposition of $T_H=W_1\cup\ldots\cup W_{k+1}$ provided by Lemma \ref{cut}. If $\eta_{T_H}\neq 0$, then each of the $W_i$ is never grafted to the last leaf of another $W_j$.
\end{lem}
\begin{proof}
If some $W_i$ were grafted to the last leaf of another $W_j$, it would lead to a composition either of the form
\begin{align}
 (\Pi_A\otimes H_B)\circ \mu^\otimes_{2,0}(\rho \otimes \tau, \Pi_A(\varphi)\otimes H_B(\psi))\nonumber
\end{align}
or $(P_A\otimes P_B)\circ \mu^\otimes_{2,0}(\rho \otimes \tau, \Pi_A(\varphi)\otimes H_B(\psi))$, depending on whether or not the root of $W_j$ coincides with the root of $T$. In either case we get zero by Lemma \ref{vanishing}(b).
\end{proof}

\begin{lem}\label{minprod}
 Let $W$ be a tree with $k$ leaves of the first type in Lemma \ref{cut} . For each $1\leq i \leq k$, consider $\xi_i=i_A(a_i)\otimes f_i$ for $a_i \in A$ and $f_i \in End_B$. If the root of $W$ does not coincide with the root $T$, then
\begin{align}
 \eta_W(\xi_1,\ldots,\xi_k)= (-1)^{\alpha + \lhd }i_A(\m_{k, \beta}(a_1,\ldots,a_k))\otimes H'_B(f_1 \circ \ldots \circ f_k),
\end{align}
where $H'(\tau)=(-1)^{\vert \tau \vert}H(\tau)$ and the sign exponent $\alpha=k + \sum_{i=1}^{k} \vert a_i \vert + (k-i-1)\vert f_i \vert$. 
If the root of $W$ does coincide with the root of $T$, then
 \begin{align}
 \eta_W(\xi_1,\ldots,\xi_k)= (-1)^{\alpha' + \lhd }\m_{k, \beta}(a_1,\ldots,a_k)\otimes f_1 \circ \ldots \circ f_k(e_B),
\end{align}
with $\alpha'=\sum_{i=1}^{k} (k-i-1)\vert f_i \vert$.

When $W$ is a tree of second type (in Lemma \ref{cut}), we have 
$$\eta_W(a\otimes b)= (-1)^{\vert a \vert}F^A_{1,0}(a)\otimes F^B_{1,\beta''}(b).$$
\end{lem}
\begin{proof}
The proof is a straightforward computation and is left to the reader.
\end{proof}

Given $T_H\in\Gamma_n(\beta_1 \times \beta_2)$, $n\geq 2$, consider the decomposition 
\begin{align}\label{decomp}
 T_H=S_1\cup\ldots\cup S_{t+1}\cup R_1\cup\ldots\cup R_j
\end{align}
provided by Lemma \ref{cut}, where each $S_i$ is of the first type and $R_l$ is of the second type. Recall that each $S_i$ has at most one vertex of valency one with weight $\beta'_i\times 0$ and each $R_l$ has a unique vertex of weight $0\times\beta''_l$. Moreover, $\sum_{i=1}^{k} \beta'_i = \beta_1$ and
$\sum_{i=1}^{k} \beta''_i = \beta_2$. 

Let $U$ be the tree obtained  from $T_H$ by forgetting  all the subtrees $R_l$ and replacing each $S_i$ by the corolla of the respective size with $\beta'_i$ assigned to its unique vertex. Then $U$ has $t$ internal edges and  $U \in L_n$ by Lemma \ref{lastleaf}. Observe that $U_{max}$ can be obtained from $T_H$ by forgetting the $R_l$ and replacing each $S_i$ by the maximal binary tree with the same number of leaves.

Applying Lemma \ref{minprod} to all subtrees $S_i$ and $R_j$ we conclude
\begin{align}
 \eta_{T_H}(a_1 \otimes b_1, \ldots , a_n \otimes b_n)= (-1)^{\dagger}U^A(a_1,\ldots,a_n)\otimes \psi(b_1,\ldots,b_n)
\end{align}
for some sign $\dagger$, which we will describe later, and for a map $\psi$ which we now describe. 

We first consider the case when there is a single $S$. Then by Lemma \ref{lastleaf}, there can be at most $n-1$ trees of second type $R_1,\ldots,R_{n-1}$ grafted to the first $n-1$ leaves of $S$. An easy computation using Lemma \ref{vanishing}(d) shows that 
\begin{align} \label{maxtree}
 \psi(b_1,\ldots,b_n)=\m_{2,\beta''_1}(b_1,\m_{2,\beta''_2}(b_2,\ldots, \m_{2,\beta''_{n-1}}(b_{n-1},b_n)\ldots)).
\end{align}
If there are several $S_i$'s, we order them so that $i<j$ whenever $S_i$ is above $S_j$. Then the root of $S_1$ agrees with the root of $T$. If we denote  the number of leaves of $S_1$ by $l$, then by definition there are $\varphi_m, \tau_k \in End_B$ so that
$$\psi(b_1,\ldots,b_n)=(\varphi_1\circ\ldots\circ\varphi_l)_0(e_B).$$
Lemma \ref{minprod} implies that for each $i$, $\varphi_i=H'(\tau_1\circ\ldots\circ \tau_{i_0})$ or $\varphi_i=F^B_{1,\beta_{i_0}}(b_{i_0})$; moreover, we must have $\varphi_l=F^B_{1,0}(b_n)$, by Lemma \ref{lastleaf}.

Let $j$ be the minimum such that $\varphi_i=F^B_{1,\beta_{i_0}}(b_{i_0})$ for $i>j$. Then, proceeding as in the calculation of (\ref{maxtree}) we have
\begin{align}
\psi (& b_1,\ldots ,b_n)=\nonumber\\
&=[\varphi_1\circ\ldots\circ\varphi_{j-1}\circ H'(\tau_1,\ldots\tau_{j_0})](\m_{2,\beta_q}(b_q,\m_{2,\beta_{q+1}}(b_{q+1},\ldots \m_{2,\beta_{n-1}}(b_{n-1},b_n))\ldots)\nonumber\\
&=\varphi_1 \circ \ldots\circ \varphi_{j-1} \Big( \tau_1 \circ \ldots \circ \tau_{j}\big( b_{q-1},\m_{2,\beta_q}(b_q,\ldots \m_{2,\beta_{n-1}}(b_{n-1},b_n)\ldots)\big)\Big),\nonumber
\end{align}
since $\tau_{j_0}=F^B_{1,0}(b_{q-1})$ by Lemma \ref{lastleaf}. 

Repeating this argument until we exhaust all of the $S_i$ and using the formula
\begin{align}
(\varphi_1\circ\ldots\circ\varphi_l)_r&(v;y_1,\ldots,y_r)=\nonumber\\
&\sum_{0\leq i_l\leq \ldots\leq i_1=r}\varphi_1(\ldots\varphi_{l_1}(\varphi_l(v; y_1,\ldots,y_{i_l});\ldots y_{i_{l-1}})\ldots),\label{iterate}
\end{align}
we conclude
\begin{align}\label{soma*}
\psi(b_1,\ldots,b_n)=\sum F^B_{1,\beta''_1}(b_1)(\ldots F^B_{1,\beta''_l}(b_l)(b_{l+1}\ldots b_q)\ldots).
\end{align}
The terms in the sum are in one to one correspondence with binary trees $Z$ with $t$ right leaning internal edges and $U_{max}\leq Z$. This can be seen as follows: Each of the $t$ right leaning internal edges in $U_{max}$ corresponds to one of the internal edges of $T_H$ labeled with $\Pi_A \otimes H_B$. Therefore it corresponds to applying $H'$ when computing $\psi$. Each time we apply $H'$, we sum over all the ways of composing elements in $End_B$ using (\ref{iterate}). These are in one to one correspondence with all the trees $W$ that can be obtained from $U_{max}$ by performing a sequence of moves described in Figure \ref{figordertree}, as long as the edge being moved is an internal edge, which are exactly the ones with $t$ right leaning internal edges and $U_{max}\leq Z$. For example, $U_{max}$ corresponds to taking $0=i_l=\ldots=i_2, \ i_1=r$ in (\ref{iterate}), every time. 

Given such tree $Z$, we collapse all the right leaning edges and obtain a tree $Q$. For each vertex $v$ in $Q$, let $p_v$ be the left most leaf under $v$. Define the weight of $v$ to be $\beta''_l$ if $p_v$ is contained in $R_l$ (for some $l$) and zero otherwise. Then the corresponding term in the sum (\ref{soma*}) is exactly $Q^B(b_1,\ldots, b_n)$. 

Thus we conclude that 
\begin{align}
 \eta_{T_H}(a_1 \otimes b_1, \ldots , a_n \otimes b_n)= (-1)^{\dagger}U^A(a_1,\ldots,a_n)\otimes \alpha(U)^B(b_1,\ldots,b_n).
\end{align}

The only thing left to show is that $\dagger$ coincides with the sign described in (\ref{tensor}). From now on all the computations are modulo 2. 

Each of the trees $S_i$ contributes to $\dagger$ with the sign described in Lemma \ref{minprod}. By Lemma \ref{lastleaf}, the last leaf of each $S_i$ coincides with one leaf of $T$, which we denote by $m_i$. This contributes with the sign $(-1)^{\vert b_{m_i} \vert}$, coming from the computation $F^B_{1,0}(b_{m_i})(e_B)=(-1)^{\vert b_{m_i} \vert}b_{m_i}$.
Therefore, noting that each $S_i$ corresponds to one vertex in $U$, we have
\begin{align}
 \dagger=&\sum_{v \in V(U)} \Big[ val(v) + \sum_{i<j} \vert \tilde{f_i} \vert \vert \tilde{a_j} \vert + \sum_{i=1}^{val(v)} \vert \tilde{a_i}\vert + (val(v)+1-i)\vert \tilde{f_i}\vert\Big]\nonumber \\
 &+ (\sum_{i=1}^{n} \vert a_i \vert + n + \vert E_{int}(U) \vert )+ \sum_{i=1}^{t+1} \vert b_{m_i}\vert,
\end{align}
where 
\begin{align}
 \vert \tilde{f_i} \vert = \sum_{p_{v(i)}\leq s \leq q_{v(i)}} \vert b_s \vert + I_{v(i)} \ \ \ \textrm{and} \ \  \vert \tilde{a_i} \vert= \sum_{p_{v(i)}\leq s \leq q_{v(i)}}  \vert a_s \vert + \tilde{I}_{v(i)} + J_{v(i)}.\nonumber
\end{align}

We will divide $\dagger= \mathcal{A} + \mathcal{B} + \mathcal{R}$ into three sums, the first involving terms with $\vert a_m\vert$, the second involving terms with $\vert b_m\vert$ and the third one with the remaining terms. Explicitly we have
\begin{align}
 \mathcal{A}=& \sum_{v \in V(U)} \Big[ \sum_{i<j} I_{v(i)}\big( \sum_{p_{v(j)}\leq s \leq q_{v(j)}} \vert a_s\vert \big) + \sum_{p_{v}\leq s \leq q_{v}} \vert a_s\vert \Big] + \sum_{i=1}^{n} \vert a_i \vert \nonumber
\end{align}

\begin{align}
 \mathcal{B}= \sum_{v \in V(U)} \Big[ \sum_{i<j} \vert b_i \vert (\tilde{I}_{v(j)}+J_{v(j)}) +\sum_i (val(v)-1-i)\sum_{p_{v(i)}\leq s \leq q_{v(i)}} \vert b_s \vert \Big] + \sum_{i=1}^{k} \vert b_{m_i}\vert \nonumber
\end{align}

\begin{align}
 \mathcal{R}=&\sum_{v \in V(U)} \big( \sum_{i} ((val(v)-i-1)I_{v(i)} + \tilde{I}_{v(i)} +J_{v(i)}) + \sum_{i<j} I_{v(i)}(\tilde{I}_{v(j)} + J_{v(j)})\big) \nonumber\\
   &+ \sum_{v \in V(U)} val(v) + n + \vert E_{int}(U)\vert \nonumber\\
   &+ \sum_{v \in V(U)} \sum_{i< j} \Big( \sum_{p_{v(i)}\leq s \leq q_{v(i)}} \vert b_s \vert\Big) \Big( \sum_{p_{v(i)}\leq r \leq q_{v(i)}} \vert a_r \vert\Big). \nonumber
  \end{align}
First observe that $\mathcal{R}=\theta(U) +\lhd$, because $\sum_{v \in V(U)} val(v) + n + \vert E_{int}(U)\vert=0$. Next we claim that $\mathcal{A}= \vert {\bf a}\vert \vert E_{int}(U) \vert + \gamma_{\bf a}$. This claim is equivalent to 
\begin{align}
 \sum_{v \in V(U)} \Big[ \sum_{i<j} I_{v(i)}\big( \sum_{p_{v(j)}\leq s \leq q_{v(j)}} \vert a_s\vert \big) + \sum_{ s \leq q_{v}} \vert a_s\vert \Big] + (\vert E_{int}(U) \vert + 1)\vert {\bf a} \vert = 0. \nonumber
\end{align}
Noting that $U$ has $\vert E_{int}(U) \vert + 1$ internal vertices, this is in turn equivalent to 
\begin{align}
 \sum_{v \in V(U)} \sum_{ s> q_{v}} \vert a_s\vert=\sum_{v \in V(T)} \sum_{i<j} I_{v(i)}\big( \sum_{p_{v(j)}\leq s \leq q_{v(j)}} \vert a_s\vert \big).\nonumber
\end{align}
Now this statement follows from the fact that both sides equal $$\sum_{i=1}^n \vert a_i\vert \nu_i,$$
where $\nu_i$ is the number of internal vertices (or internal edges) of $U$ to the left of the i-th leaf.

At last we show that $\mathcal{B}=\gamma_{{\bf b}}$. As above we can give an alternative description of $\gamma_{{\bf b}}$, namely
$$\gamma_{{\bf b}}=\sum_{i=1}^n \vert b_i\vert \bar{\nu}_i,$$
where $\bar{\nu}_i$ is the number of internal vertices in $\alpha_0(U)$ to the right of the i-th leaf. If we exclude the internal vertex adjacent to the root (which is not to the right of any leaf), internal vertices in $\alpha_0(U)$ are in one-to-one correspondence with the left leaning edges in $U_{max}$. These correspond to the internal edges of the $S_m$ in (\ref{decomp}).
On the other hand we can rewrite $\mathcal{B}$ as
$$\mathcal{B}=\sum_{v \in V(T)} \sum_{i<j} \vert b_i \vert (\tilde{I}_{v(j)}+J_{v(j)}) +\sum_{i}^{val(v)-1} (val(v)-1-i)\sum_{p_{v(i)}\leq s \leq q_{v(i)}} \vert b_s \vert. $$
If we denote by $S_v$ the subtree in $T_H$ corresponding to $v \in V(U)$, then $(val(v)-1-i)$ is the number of internal edges of $S_v$ to the right of its $i^{th}$ leaf. Also $(\tilde{I}_{v(j)}+J_{v(j)})$  is the number of internal edges of the other $S_{v'}$, under $v(j)$. Therefore we conclude that 
$$\mathcal{B}=\sum_{i=1}^n \vert b_i\vert \bar{\nu}_i=\gamma_{{\bf b}}.$$
This completes the proof of Theorem \ref{filteredtensor}. 
\end{proof}

\section{Criterion for tensor product}

In this section we will describe a set of conditions under which a given filtered \Ai-algebra is the tensor product of two subalgebras. Although rather restrictive, this criterion will be sufficient for the two applications we have in mind, namely, the description of deformations of a tensor product via  bounding cochains in Section 6 and the proof of the K\"unneth Theorem for the Fukaya algebra of the product of two Lagrangian submanifolds in \cite{Amo}.
%This criterion is inspired by the moduli spaces of holomorphic disks that were used in \cite{Amothesis}. These moduli spaces are essentially the moduli spaces of quilted strips.

We start by defining the notions of subalgebra and commuting subalgebras.
\begin{defn}\label{subalgebra}
Let $(A,\m^A)$ and $(C,\mu)$ be (respectively) $G_A$ and $G$-gapped filtered \Ai-algebras, for discrete submonoids $G_A\subseteq G$. We say A is a subalgebra of $C$ if $A\subseteq C$, $e_A=e_C$ and for all $k>0$ and $a_1, \ldots a_k \in A$ we have
\begin{align}
 \mu_{k,\beta}(a_1,\ldots,a_k)&=\m_{k,\beta}^A(a_1,\ldots,a_k), \ \  \beta\in G_A,\nonumber \\
 \mu_{k,\beta}(a_1,\ldots,a_k)&=0, \  \ \beta\in G\setminus G_A.\nonumber
\end{align}

\end{defn}
\begin{defn}\label{comsubalg}
Let $(A,\m^A)$ and $(B,\m^B)$ be $G_A$ and $G_B$-gapped filtered \Ai-algebras. Suppose $A$ and $B$ are subalgebras of $(C,\mu)$ a $G$-gapped \Ai-algebra with $G=G_A+G_B$. Define a map $K:A\otimes B\lto C$ by $K(a\otimes b)=(-1)^{\vert a\vert}\mu_{2,0}(a,b)$. We say $A$ and $B$ are \emph{commuting subalgebras} if given $c=K(a\otimes b)$ and $c_1,\ldots, c_k\in C$ such that for each $i$, $c_i=a_i$ or $c_i=b_i$ for some $a_i\in A$ and $b_i\in B$, the following conditions hold.
\begin{itemize}
\setlength{\parsep}{0pt}
\setlength{\itemsep}{0pt}
\item[{\bf(a)}] For $k>0$, $\mu_{k,\beta}(c_1,\ldots, c_k)=0$ unless
\begin{itemize}
\setlength{\parsep}{0pt}
\setlength{\itemsep}{0pt}
\item[{\bf(i)}] $(k,\beta)=(2,0)$ and $c_1 \in A$, $c_2 \in B$ (or vice-versa) in which case, $\mu_{2,0}(c_1,c_2)+(-1)^{\|c_1\|\|c_2\|}\mu_{2,0}(c_2,c_1)=0$,
\item[{\bf(ii)}] $c_i=a_i$ for all $i$ and $\beta\in G_1$,
\item[{\bf(iii)}] $c_i=b_i$ for all $i$ and $\beta\in G_2$.
\end{itemize}
\item[{\bf(b)}] $\mu_{0,\beta}=\m^A_{0,\beta}+\m^B_{0,\beta}$, with the convention that $\m^A_{0,\beta}=0$ (respectively $\m^B_{0,\beta}$) if $\beta\notin G_A$ (respectively $\beta\notin G_B$).
\item[{\bf(c)}] $\mu_{k+1,\beta}(c_1,\ldots,c_i,c,c_{i+1},\ldots,c_k)=0$ unless
\begin{itemize}
\setlength{\parsep}{0pt}
\setlength{\itemsep}{0pt}
\item[{\bf(i)}] $c_i=a_i$  for all $i$ and $\beta\in G_A$, in which case it equals
$$(-1)^{|b|\sum_{j> i}\|a_j\|}K(\m^A_{k+1,\beta}(a_1\ldots a_i,a,\ldots,a_k)\otimes b),$$
\item[{\bf(ii)}] $c_i=b_i$ for all $i$ and $\beta\in G_B$, in which case it equals
$$(-1)^{|a|\left(\sum_{j\leq i}\|b_j\|+1\right)}K(a\otimes \m^B_{k+1,\beta}(b_1\ldots b_i,b,\ldots,b_k)).$$
\end{itemize}
\end{itemize}
\end{defn}
We are now ready to state the main theorem of this section.
\begin{thm}\label{criterion}
Let $(A,\m^A)$ and $(B,\m^B)$ be commuting subalgebras of $(C,\mu)$ and equip $A\otimes B$ with the differential $\delta=\m_{1,0}^A\otimes id+id\otimes \m_{1,0}^B$. If $K:A\otimes B\lto C$ is an injective map that induces an isomorphism on $\delta - \mu_{1,0}$-cohomology then there is a (strict) quasi-isomorphism
$$A\otimes_\infty B \simeq C.$$ 
\end{thm}
We will prove this theorem in three steps. First we will replace $(C,\mu)$ by a quasi-isomorphic \Ai-algebra $(A\otimes B, \eta)$ whose underlying vector space is $A\otimes B$. Second, we will construct a new filtered dg-algebra $End_{A,B}$ quasi-isomorphic to $A\otimes_\infty B$ and finally, we will construct a quasi-isomorphism from $(A\otimes B, \eta)$ to $End_{A,B}$.
\begin{prop}\label{restrict}
Let $(A,\m^A)$, $(B,\m^B)$ and $(C,\mu)$ be as in Theorem \ref{criterion}. Then there is a filtered \Ai-algebra $(A\otimes B, \eta)$ quasi-isomorphic to $(C,\mu)$, such that $A$ and $B$ are commuting subalgebras of $(A\otimes B, \eta)$, via the inclusions $a\lmto a\otimes e_B$ and $b\lmto e_A\otimes b$.
\end{prop}
\begin{proof}
The proof is an application of the homological perturbation lemma. By assumption on the map $K$, we can choose a subspace $V\subset C$ such that $$C=K(A\otimes B)\oplus V.$$
Using this decomposition we write $$\mu_{1,0}=\left(\begin{array}{cc}\delta&g\\0&d_V
\end{array}\right).$$
Since $K$ is an isomorphism on cohomology, $d_V$ must be acyclic, thus we can choose $W$ such that
$$V=W\oplus d_V(W).$$
We define $p:C\lto A\otimes B$ and $H:C\lto C$ by
\begin{align}
p(K(a\otimes b),w_0,d_Vw_1)&=a\otimes b-g(w_1),\nonumber\\
H(K(a\otimes b),w_0,d_Vw_1)&=(0,-w_1,0).\nonumber
\end{align}
We can easily check that $K$, $p$ and $H$ define homotopy data satisfying the side conditions. Therefore, applying Theorem \ref{hpl} and Proposition \ref{sideconditions} we obtain a filtered \Ai-algebra $(A\otimes B, \eta)$ quasi-isomorphic to $(C,\mu)$.

Next we show that $A$ is a subalgebra of $(A\otimes B, \eta)$. Recall that 
$$\eta_{k,\beta}(a_1\otimes e, \ldots, a_k\otimes e)=\sum_{T\in\Gamma_k(\beta)}\eta_{T}(a_1\otimes e, \ldots, a_k\otimes e).$$
Since by assumption $A\subseteq A\otimes B\subseteq C$ is closed under the $\mu_{k,\beta}$ operations and $H|_{A\otimes B}=0$, we see that if $T$ has internal edges, then $\eta_{T}(a_1\otimes e, \ldots, a_k\otimes e)=0$. Since $A$ is a subalgebra of $C$, this together with $K(a\otimes e)=a$ implies $\eta_{k,\beta}(a_1\otimes e, \ldots, a_k\otimes e)=\mu_{k,\beta}(a_1,\ldots,a_k)=\m_{k,\beta}^A(a_1,\ldots,a_k)$. Similarly, we prove that $B$ is  a subalgebra of $(A\otimes B, \eta)$.

Next we verify that the conditions in Definition \ref{comsubalg} are satisfied. The same argument implies that
$$\eta_{k,\beta}(c_1,\ldots,c_k)=\mu_{k,\beta}(c_1,\ldots,c_k),$$
$$\eta_{k+1,\beta}(c_1,\ldots,c_i,c,c_{i+1},\ldots,c_k)=\mu_{k,\beta}(c_1,\ldots,c_i,c,c_{i+1},\ldots,c_k).$$
It now easily follows that $(A,\m^A)$ and $(B,\m^B)$ are commuting subalgebras of $(A\otimes B, \eta)$.
\end{proof}

Next we introduce a new filtered dg-algebra $End_{A,B}$. As a vector space, $End_{A,B}$ is the subspace of
$$\Hom_\mathbb{K}\left(\bigoplus_{r,s\geq 0}A\otimes B\otimes A^{\otimes r}\otimes B^{\otimes s}, A\otimes B\right)$$
satisfying
$$\rho_{r,s}(\bullet;\ldots,e_A,\ldots)=\rho_{r,s}(\bullet;\ldots,e_B,\ldots)=0.$$
An element $\rho=\{\rho_{r,s}\}_{r,s}\in End_{A,B}$ is said to have degree $k=|\rho|$ if each $\rho_{r,s}$ has degree $k-r-s$.

To define the \Ai \ operations on $End_{A,B}$, we introduce the following convenient notation
\begin{align}\label{barmk}
\overline \m^A_{k+1,\beta}(a\otimes b,a_1,\ldots,a_k)&=(-1)^{|b|\sum_{i=1}^{k}||a_i||}\m^A_{k+1}(a,a_1,\ldots,a_k)\otimes b,\nonumber\\
\overline \m^B_{k+1,\beta}(a\otimes b,b_1,\ldots,b_k)&=(-1)^{|a|}a\otimes \m^B_{k+1}(b,b_1,\ldots,b_k).
\end{align} 
For each $\beta\in G=G_A + G_B$, we define
$$\overline \mu _{k,\beta}=\sum_{\substack{\beta_1\in G_A, \beta_2 \in G_B\\ \beta_1+\beta_2=\beta}}\overline \mu _{k,\beta_1\times\beta_2},$$ 
where the maps $\overline\mu _{k,\beta_1\times\beta_2}$ are defined  by
\begin{align*}
(\overline\mu _{0,\beta_1\times 0})_{r,s}(x\otimes y;\ba;\bb)&=\left\{
\begin{array}{cl}
\sum_{\beta_1'+\beta_1''=\beta_1}(-1)^{|y|||\ba||}\m^A_{r+2,\beta_1'}(\m^A_{0,\beta_1''},x,\ba)\otimes y& \  s=0,\\
0&\textrm{otherwise},
\end{array}\right.\\
(\overline\mu _{0,0\times\beta_2})_{r,s}(x\otimes y;\ba;\bb)&=\left\{
\begin{array}{cl}
\sum_{\beta_2'+\beta_2''=\beta_2}x\otimes \m^B_{s+2,\beta_2'}(\m^B_{0,\beta_2''},y,\bb)& \  r=0,\\
0&\textrm{otherwise},
\end{array}\right.
\end{align*}
\begin{align*}
\overline\mu _{1,\beta_1\times 0}(\rho)_{r,s}(& x\otimes y; \ba;\bb)=\sum(-1)^{1+||\ba^{(2)}|| ||\bb||}\overline \m^A_{r-i+1,\beta_1}\left(\rho_{i,s}(x\otimes y;\ba^{(1)};\bb),\ba^{(2)}\right)\\
&+\sum(-1)^{|\rho|}\rho_{r-i,s}\left(\overline \m^A_{i+1,\beta_1}(x\otimes y;\ba^{(1)});\ba^{(2)};\bb\right)\\
&+\sum(-1)^{|\rho|+||x\otimes y||+||\ba^{(1)}||}\rho_{r-j+1,s}\left(x\otimes y;\ba^{(1)},\m^A_{j,\beta_1}(\ba^{(2)}),\ba^{(3)};\bb\right),\\
\overline\mu _{1,0\times\beta_2}(\rho)_{r,s}(& x\otimes y; \ba;\bb)=\sum-\overline \m^B_{j+1,\beta_2}\left(\rho_{r,s-j}(x\otimes y;\ba;\bb^{(1)}),\bb^{(2)}\right)\\
&+\sum(-1)^{|\rho|+||\ba||\ ||\bb^{(1)}||}\rho_{r,s-j}\left(\overline \m^B_{j+1,\beta_2}(x\otimes y;\bb^{(1)});\ba;\bb^{(2)}\right)\\
&+\sum(-1)^{|\rho|+||x\otimes y||+||\ba||+||\bb^{(1)}||}\rho_{r,s-j+1}\left(x\otimes y;\ba;\bb^{(1)},\m^B_{j,\beta_2}(\bb^{(2)}),\bb^{(3)}\right).
\end{align*}
$$\overline \mu _{2,0}(\rho,\tau)(x\otimes y;\ba;\bb)=\sum(-1)^{|\rho|+||\bb^{(1)}||\ ||\ba^{(2)}||}\rho_{r-i,s-j}\left(\tau_{i,j}\left(x\otimes y;\ba^{(1)};\bb^{(1)}\right); \ba^{(2)};\bb ^{(2)}\right)$$ 
and all other $\overline \mu_{k,\beta_1\times\beta_2}=0$. 
We have the following:
\begin{prop}
$End_{A,B}=\left(End_{A,B},\overline \mu \right)$ is a filtered dg-algebra.
\end{prop}
\begin{proof}
The proof is straightforward computation combining the proofs of Lemmas \ref{EndAdga} and \ref{tensordg}. We simply highlight the main points.

Just as in Lemma \ref{EndAdga} we can show that
$$\sum_{\beta_1'+\beta_1''=\beta_1}\overline\mu _{1,\beta_1'\times 0}\left(\overline\mu _{0,\beta_1''\times 0}\right)=0=\sum_{\beta_2'+\beta_2''=\beta_2}\overline\mu _{1,0\times\beta_2'}\left(\overline\mu _{0,0\times\beta_2''}\right).$$
This, combined with the equalities
$$\overline\mu _{1,\beta_1\times 0}\left(\overline\mu _{0,0\times\beta_2}\right)=\overline\mu _{1,0\times\beta_2}\left(\overline\mu _{0,\beta_1\times 0}\right)=0,$$
immediately implies the \Ai-equation with no inputs. For the next equation, we note that
$$\overline\mu _{1,\beta_1\times 0}\left(\overline\mu _{1,0\times\beta_2}(\rho)\right)+\overline\mu _{1,0\times\beta_2}\left(\overline\mu _{1,\beta_1\times 0}(\rho)\right)=0.$$
As in Lemma \ref{EndAdga}, we can verify the equation
$$\sum_{\beta_1'+\beta_1''=\beta_1}\overline\mu _{1,\beta_1'\times 0}\left(\overline\mu _{1,\beta_1''\times 0}(\rho)\right)+\overline\mu _{2,0}\left(\overline\mu _{0,\beta_1\times 0},\rho\right)+(-1)^{||\rho||}\overline\mu _{2,0}\left(\rho,\overline\mu _{0,\beta_1\times 0}\right)=0$$
and the analog for $0\times\beta_2$. These combine to prove the \Ai \ equation with one input. Finally, we can prove
$$\overline\mu_{1,\gamma}\left(\overline\mu _{2,0}(\rho,\tau)\right)+\overline\mu _{2,0}\left(\overline\mu _{1,\gamma}(\rho),\tau\right)+(-1)^{||\rho||}\overline\mu _{2,0}\left(\rho,\overline\mu _{1,\gamma}(\tau)\right)=0$$ for $\gamma=\beta_1\times 0$ or $\gamma=0\times\beta_2$, and
$$\overline\mu_{2,0}\left(\overline\mu _{2,0}(\rho,\tau),\eta\right)+(-1)^{||\rho||}\overline\mu_{2,0}\left(\rho,\overline\mu _{2,0}(\tau,\eta)\right)=0,$$ by straightforward computations. This finishes the proof of the \Ai \ equations. To complete the proof, we note that
\begin{align*}
(id)_{r,s}(x\otimes y;\ba;\bb)=\left\{
\begin{array}{cl}
x\otimes y&,r=s=0,\\
0&,\textrm{ otherwise},
\end{array}\right.
\end{align*}
is a unit for $End_{A,B}$.
\end{proof}

\begin{prop}\label{EndAB}
Let $A$ and $B$ be filtered \Ai-algebras. The map 
$$S(\rho\otimes \tau)(x\otimes y ;\ba;\bb)=(-1)^{|y|||\ba||+|\tau|(|x|+||\ba||)}\rho(x;\ba)\otimes\tau(y;\bb)$$ 
defines a naive quasi-isomorphism $S:End_A\otimes_{dg}End_B\lto End_{A,B}$.
\end{prop}
\begin{proof}
It is straightforward check that
$$S(\mu_{0,\gamma}^\otimes)=\overline{\mu}_{0,\gamma},$$
$$S(\mu_{1,\gamma}^\otimes(\varphi))=\overline{\mu}_{1,\gamma}(S(\varphi))$$
and
$$S(\mu_{2,\gamma}^\otimes(\varphi_1,\varphi_2))=\overline{\mu}_{2,0}(S(\varphi_1),S(\varphi_2)),$$
for any $\varphi=\rho\otimes\tau$ and $\gamma=\beta_1\times 0$ or $\gamma=0\times\beta_2$.
This shows that $S$ is a naive \Ai-homomorphism. We need to check that the induced map
$$S_*:H^*(End_A\otimes_{dg}End_B,\mu_{1,0}^\otimes)\lto H^*(End_{A,B},\overline{\mu}_{1,0})$$
is an isomorphism. In fact, when $A$ and $B$ are finite dimensional vector spaces $S$ is an isomorphism. In the general case we argue as follows. Let $I$ denote the composition
$$A\otimes B\xrightarrow{F^A_{1,0}\otimes F^B_{1,0}}End_A\otimes End_B\stackrel{S}{\lto}End_{A,B}.$$
Lemma \ref{homotopytensor} states that $F^A_{1,0}$ and $F^B_{1,0}$ induce isomorphisms in cohomology. Therefore $F^A_{1,0}\otimes F^B_{1,0}$ also induces an isomorphism in cohomology by the K\"unneth theorem. We have reduced the proof to showing that $I$ induces an isomorphism in cohomology. We will do this by providing a homotopy inverse. Define $\bar P:End_{A,B}\lto A\otimes B$ by
$$\bar P(\rho)=(-1)^{|\rho|}\rho_{0,0}(e_A\otimes e_B).$$ 
By definition,
$$I(u\otimes v)(x\otimes y;\ba;\bb)=(-1)^{|y|||\ba||+(|x|+||a||)|v|}\m^A_{r+2,0}(u,x,\ba)\otimes \m^B_{s+2,0}(v,y,\bb),$$
which readily implies $\bar P\circ I=id_{A\otimes B}$. The composition $I\circ \bar P$ is homotopic to the identity in $End_{A,B}$. In fact, if we define
\begin{align*}
(K\rho)_{r,s}(x\otimes y;\ba;\bb)&=(-1)^{\xi_1}\overline \m^A_{r+2,0}\left(\rho_{0,s+1}(E_a\otimes e_B;y;\bb),x,\ba\right)\\
&+(-1)^{\xi_2}\rho_{r+1,s}(e_A\otimes y;x,\ba;\bb),
\end{align*}
where $\xi_1=||\bb||(||\ba||+||x||)+|x||y|+||y||$ and $\xi_2=|\rho|+|y|||x||$, we can see that
$$K\overline{\mu}_{1,0}+\overline{\mu}_{1,0}K=I\circ\bar P-id_{End_{A,B}}.$$ 
This is a simple albeit long computation, entirely analogous to Lemma \ref{homotopytensor}, so we will omit it.
\end{proof}

The third and main step in the proof of Theorem \ref{criterion} is the construction of an \Ai-homomorphism
$$F:(A\otimes B,\eta)\lto End_{A,B}.$$
For this, we need to introduce the notion of shuffle product.
\begin{defn}
Let $r,s\geq 0$ be integers. We say a permutation $\sigma\in S_{r+s}$ is a \emph{$(r,s)$-shuffle} if $\sigma(i)<\sigma(j)$ for $i<j\leq r$ or $r+1\leq i<j$.
\end{defn}
Consider $\ba=a_1\otimes\ldots\otimes a_r\in A^{\otimes r}$ and $\bb=b_1\otimes\ldots\otimes b_s\in B^{\otimes s}$ and denote
\begin{align*}
u_i=\left\{\begin{array}{ll}
a_i&,\  1\leq i\leq r,\\
b_{i-r}&,\ r<i\leq r+s.
\end{array}\right.
\end{align*}
Given a $(r,s)$-shuffle $\sigma$, we define
$$\sigma(\ba,\bb)=(-1)^\epsilon u_{\sigma^{-1}(1)}\otimes\ldots\otimes u_{\sigma^{-1}(r+s)},$$ with $$\epsilon =\sum_{\substack{i\leq r<j\\ \sigma(j)<\sigma(i)}}||u_j||||u_i||.$$With the same notation, we define
$$Sh(\ba;\bb)=\sum_{\sigma}\sigma(\ba;\bb)$$ where we sum over all $(r,s)$-shuffles $\sigma$.
In the next lemma we establish some properties of the shuffle $Sh$ that we will need later.

\begin{lem}\label{shprop}
Let $A$ and $B$ be commuting subalgebras of $(A\otimes B,\eta_{k,\beta})$ as in Proposition \ref{restrict}. We have the following:
\begin{itemize}
\item[{\bf(a)}]   Let $\alpha=||\ba^{(2)}||||\bb^{(1)}||$ we have 
$$\sum_{(Sh(\ba;\bb))} Sh(\ba;\bb)^{(1)}\otimes Sh(\ba;\bb)^{(2)}=\sum_{(\ba)(\bb)}(-1)^{\alpha}Sh(\ba^{(1)};\bb^{(1)})\otimes Sh(\ba^{(2)};\bb^{(2)}).$$

\item[{\bf(b)}] $ Sh(\ba_1,x,\ba_2;\bb)=\sum_{(\bb)}(-1)^{||\bb^{(1)}||(||x||+||\ba^{(2)}||)} Sh(\ba_1;\bb^{(1)})\otimes x \otimes Sh(\ba_2;\bb^{(2)})$ 
and $Sh(\ba;\bb_1,y,\bb_2)=\sum_{(\ba)}(-1)^{||\ba^{(2)}||(||\bb^{(1)}||+||y||)} Sh(\ba^{(1)};\bb_1)\otimes y \otimes Sh(\ba^{(2)};\bb_2).$
\item[{\bf(c)}] $\eta_{r+s+1,\beta}(z,Sh(\ba;\bb))=\left\{
\begin{array}{cl}\overline \m^A_{r+1,\beta\times 0}(z,\ba)& \ \beta\in G_A,\ s=0,\\
			\overline \m^B_{s+1,0\times \beta}(z,\bb)& \ \beta\in G_B,\ r=0,\\
			\overline \m^A_{1,\beta\times 0}(z)+\overline \m^B_{1,0\times \beta}(z)& \ r,s=0,\\
			0& \ \textrm{otherwise.}
\end{array}\right.$

\item[{\bf(d)}] $\eta_{r+s,\beta}(Sh(\ba;\bb))=\left\{
\begin{array}{cl} \m^A_{r,\beta}(\ba)& \ \beta\in G_A,\ s=0,\ r>0,\\
			 \m^B_{s,\beta}(\bb)& \ \beta\in G_B,\ r=0,\ s>0,\\
			 \m^A_{0,\beta}(z)+\m^B_{0,\beta}(z)& \ r,s=0,\\
			 0& \ \textrm{otherwise.}
\end{array}\right.$

\end{itemize}
%with the convention $\m^A_{0,\beta}=0$ (respectively $\m^B_{0,\beta}$) if $\beta\notin G_A$ (respectively $\beta\notin G_B$).
\end{lem}
\begin{proof}
Part (a) follows from the fact that for any vector space $V$, the tensor algebra $TV=\oplus_{k\geq0}V^{\otimes k}$, equipped with $\Delta(v)=\sum_{(v)}v^{(1)}\otimes v^{(2)}$ and $Sh$ is a bialgebra. See \cite{GetJon} for a proof of this fact.

Part (b) is a straightforward check that we will omit. For (c), by definition we have
$$\eta_{r+s+1,\beta}(z,Sh(\ba;\bb))=\sum_{\sigma}\eta_{r+s+1,\beta}(z,\sigma(\ba;\bb)).$$
Since $A$ and $B$ are commuting subalgebras,
$$\eta_{r+s+1,\beta}(z,\sigma(\ba;\bb))=0$$
unless one of the following happens: $s=0$ and $\beta\in G_A$, in which case $\sigma=id$ and $\eta_{r+1,\beta}(z,\ba)=\overline \m^A_{r+1,\beta\times 0}(z,\ba)$; or $r=0$ and $\beta\in G_B$ which implies $\eta_{r+1,\beta}(z,\bb)=\overline \m^B_{s+1,0 \times \beta}(z,\bb)$. When $r=s=0$ we sum the two contributions and get $\m^A_{1,\beta\times 0}(z)+\overline \m^B_{1,0\times \beta}(z)$. 

Part (d) is the same as (c) with one extra case: $\ba=a_1$, $\bb=b_1$ and $\beta=0$. In this case
$$\eta_{2,0}(Sh(\ba;\bb))=\eta_{2,0}(a_1,b_1)+(-1)^{||b_1|| ||a_1||}\eta_{2,0}(b_1,a_1)=0,$$
since $A$ and $B$ commute. 
\end{proof}

\begin{prop}\label{mapsh}
Let $A$ and $B$ be commuting subalgebras of $(A\otimes B,\eta)$, as in Proposition \ref{restrict}. Consider the sequence of maps $F_{k,\beta}:(A\otimes B)^{\otimes k}\lto End_{A,B}$ defined by
$$F_{k,\beta}(c_1,\ldots,c_k)(z;\ba;\bb)=\eta_{k+r+s+1,\beta}(c_1,\ldots,c_k,z,Sh(\ba;\bb))$$ for $c_i, z=x\otimes y\in A\otimes B$ and $k>0$. These define a strict, filtered \Ai-homomorphism.
\end{prop}
\begin{proof}
We first check the \Ai \ equation, with no inputs:
$$\overline \mu _{0,\beta}=\sum_{\beta'+\beta''=\beta}F_{1,\beta'}(\eta_{0,\beta''}).$$
Since $A$ and $B$ are commuting subalgebras, $\eta_{0,\beta''}=\m^A_{0,\beta''}+\m^B_{0,\beta''}$, we compute
\begin{align*}&\sum_{\beta'+\beta''=\beta}F_{1,\beta'}(\eta_{0,\beta''})(z;\ba;\bb)=\sum_{\beta'+\beta''=\beta}F_{1,\beta'}(\m^A_{0,\beta''}+\m^B_{0,\beta''})(z;\ba;\bb)\\
 &=\sum_{\beta'+\beta''=\beta}\eta_{r+s+2,\beta'}\left(\m^A_{0,\beta''},z,Sh(\ba;\bb)\right)+\eta_{r+s+2,\beta'}\left(\m^B_{0,\beta''},z,Sh(\ba;\bb)\right).\end{align*}
Again, by commutativity of $A$ and $B$ this equals
$$\sum_{\beta'+\beta''=\beta}(-1)^{|y|||\ba||}\m^A_{r+2,\beta'}\left(\m^A_{0,\beta''},x,\ba\right)\otimes y+x\otimes \m^B_{s+2,\beta'}\left(\m^B_{0,\beta''},y,\bb\right) $$ with the convention that the first (respectively second) term is zero if $s\neq 0$ (respectively $r\neq 0$). This is exactly the definition of $\overline \mu _{0,\beta}(z;\ba;\bb)$.

To prove the \Ai-homomorphism equation with $k>0$ inputs we consider the \Ai-equation on $A\otimes B$, with inputs $\bc=c_1\otimes \ldots\otimes c_k$, $z$ and $\sigma (\ba;\bb)$ for some shuffle $\sigma$. Summing over $\sigma$, we obtain the equation
\begin{align}
0=& \ \sum (-1)^{*}\eta_{j',\beta'}\left(\bc^{(1)},\eta_{j'',\beta''}(\bc^{(2)}),\bc^{(3)},z,Sh(\ba;\bb)\right)+\nonumber\\
&+\sum (-1)^{*}\eta_{j',\beta'}\left(\bc^{(1)},\eta_{j'',\beta''}(\bc^{(2)},z,Sh(\ba;\bb)^{(1)}),Sh(\ba;\bb)^{(2)}\right)+\label{bigsum}\\
&+\sum (-1)^{*}\eta_{j',\beta'}\left(\bc,z,Sh(\ba;\bb)^{(1)},\eta_{j'',\beta''}(Sh(\ba;\bb)^{(2)}),Sh(\ba;\bb)^{(3)}\right).\nonumber
\end{align}
By definition, the first sum in (\ref{bigsum}) equals
\begin{align}\label{sum1}
 \sum (-1)^{||\bc^{(1)}||} F_{k-j+1,\beta'}\left(\bc^{(1)},\eta_{j,\beta''}(\bc^{(2)}),\bc^{(3)}\right)(z;\ba;\bb).
\end{align}
Applying Lemma \ref{shprop}(a) to the second sum in (\ref{bigsum}) we see that it equals 
$$\sum_{(\bc),(\ba),(\bb)} (-1)^{\zeta}\eta_{j',\beta'}\left(\bc^{(1)},\eta_{j'',\beta''}\left(\bc^{(2)},z,Sh\left(\ba^{(1)};\bb^{(1)}\right)\right),Sh\left(\ba^{(2)};\bb^{(2)}\right)\right),
$$
where $\zeta=||\bc^{(1)}||+||\bb^{(1)}||||\ba^{(2)}||$. We will further decompose this sum into two parts. The first consists of the terms where $\bc^{(1)},\bc^{(2)}\neq\bc$. It equals
\begin{align}
\sum_{\bc^{(1)},\bc^{(2)}\neq\bc} &(-1)^{\zeta}\eta_{j',\beta'}\left(\bc^{(1)},\eta_{j'',\beta''}\left(\bc^{(2)},z,Sh\left(\ba^{(1)};\bb^{(1)}\right)\right),Sh\left(\ba^{(2)};\bb^{(2)}\right)\right) \nonumber\\
 &= -\sum \overline \mu _{2,0}\left(F_{k',\beta'}\left(\bc^{(1)}\right),F_{k'',\beta''}\left(\bc^{(2)}\right)\right)(z;\ba;\bb),\label{sum2}
\end{align}
because $|F_{k,\beta}(\bc)|=||\bc||+1$. The second consists of the terms with $\bc^{(1)}=\bc$ or $\bc^{(2)}=\bc$ and is equal to
\begin{align}
&\ \sum (-1)^{||\bb^{(1)}||||\ba^{(2)}||}\eta_{j',\beta'}\left(\eta_{j'',\beta''}\left(\bc,z,Sh\left(\ba^{(1)};\bb^{(1)}\right)\right),Sh\left(\ba^{(2)};\bb^{(2)}\right)\right)\nonumber\\
&+\sum (-1)^{||\bc^{(1)}||+||\bb^{(1)}||||\ba^{(2)}||}\eta_{j',\beta'}\left(\bc,\eta_{j'',\beta''}\left(z,Sh\left(\ba^{(1)};\bb^{(1)}\right)\right),Sh\left(\ba^{(2)};\bb^{(2)}\right)\right)\nonumber\\
&=\sum_{\beta'\in G_A}(-1)^{||\bb||||\ba^{(2)}||}\overline \m^A_{j',\beta'\times 0}\left(\eta_{j'',\beta''}\left(\bc,z,Sh\left(\ba^{(1)};\bb\right)\right),\ba^{(2)}\right)\nonumber\\
&+\sum_{\beta'\in G_B}\overline \m^B_{j',0\times \beta'}\left(\eta_{j'',\beta''}\left(\bc,z,Sh\left(\ba;\bb^{(1)}\right)\right),\bb^{(2)}\right)\label{sum2b}\\
&+\sum_{\beta''\in G_A}(-1)^{||\bc||}\eta_{j',\beta'}\left(\bc,\overline \m ^A_{j'',\beta''\times 0}\left(z,\ba^{(1)}\right),Sh\left(\ba^{(2)};\bb\right)\right)\nonumber\\
&+\sum_{\beta''\in G_B}(-1)^{||\bc||+||\bb^{(1)}||||\ba||}\eta_{j',\beta'}\left(\bc, \overline \m^B_{j'',0\times \beta''}\left(z,\bb^{(1)}\right),Sh\left(\ba;\bb^{(2)}\right)\right),\nonumber
\end{align}
by Lemma \ref{shprop}(c).
Finally, we apply Lemma \ref{shprop}(a) to the third sum in (\ref{bigsum}) and get
$$\sum (-1)^\nu \eta_{j',\beta'}\left(\bc,z,Sh\left(\ba^{(1)};\bb^{(1)}\right),\eta_{j'',\beta''}\left(Sh\left(\ba^{2)};\bb^{(2)}\right)\right),Sh\left(\ba^{(3)};\bb^{(3)}\right)\right),$$ 
where $\nu=||\bc||+||z||+||Sh\left(\ba^{(1)};\bb^{(1)}\right)||+||\bb^{(1)}||\left(||\ba^{(2)}||+||\ba^{(3)}||\right)+||\bb^{(2)}||||\ba^{(3)}||$. Observe that in this expression $\eta_{j'',\beta''}\left(Sh\left(\ba^{2)};\bb^{(2)}\right)\right)$ equals $\m^A_{j'',\beta''}(\ba^{(2)})$ or $\m^B_{j'',\beta''}(\bb^{(2)})$, by Lemma \ref{shprop}(d). Then applying Lemma \ref{shprop}(b), with $x=\m^A_{j'',\beta''}(\ba^{(2)})$ and $y=\m^B_{j'',\beta''}(\bb^{(2)})$ we conclude the above expression equals
\begin{align}\label{sum3}
&\sum_{\beta''\in G_A}(-1)^{||\bc||+||z||+||\ba^{(1)}||}\eta_{j',\beta'}\left( \bc,z,Sh\left(\ba^{(1)},\m^A_{\beta''}\left(\ba^{(2)}\right),\ba^{(3)};\bb\right)\right)\nonumber\\
&+\sum_{\beta''\in G_B}(-1)^{||\bc||+||z||+||\ba||+||\bb^{(1)}||}\eta_{j',\beta'}\left( \bc,z,Sh\left(\ba,\bb^{(1)},\m^B_{\beta''}\left(\bb^{(2)}\right),\bb^{(3)}\right)\right).
\end{align}
Comparing with the definition of $\overline \mu_{1,\beta}$, we see that
\begin{align}\label{sum3b}
 (\ref{sum2b})+(\ref{sum3})=-\sum\overline \mu_{1,\beta'}\left(F_{k,\beta''}(\bc)\right)(z;\ba;\bb).
\end{align}
Assembling (\ref{sum1}), (\ref{sum2}) and (\ref{sum3b}) we conclude that (\ref{bigsum}) is equivalent to
\begin{align*}
&\sum(-1)^{||\bc^{(1)}||}F_{k-j+1,\beta'}\left(\bc^{(1)},\eta_{j,\beta''}\left(\bc^{(2)}\right),\bc^{(3)}\right)(z;\ba;\bb)\\
&-\sum \overline{\mu}_{2,0}\left(F_{k',\beta'}\left(\bc^{(1)}\right),F_{k'',\beta''}\left(\bc^{(2)}\right)\right)(z;\ba;\bb)\\
&-\sum \overline \mu _{1,\beta'}\left(F_{k,\beta''}(\bc)\right)(z;\ba;\bb)=0.
\end{align*}
This completes the proof that $F_{k,\beta}$ defines an \Ai-homomorphism. Lastly, the equality
\begin{align*}
& F_{k,\beta}\left(\bc_1,e_A\otimes e_B,\bc_2\right)(z;\ba;\bb)=\\
& =\eta_{k+r+s+1,\beta}\left(\bc_1,e_A\otimes e_B,\bc_2,z,Sh(\ba;\bb)\right)=\left\{\begin{array}{cl}z& ,(k,\beta)=(1,0),\ r=s=0,\\
0& , \textrm{otherwise},
\end{array}
\right.
\end{align*}
shows that $F_{k,\beta}$ is unital.
\end{proof}

\begin{proof}[Proof of Theorem \ref{criterion}]
Assembling Propositions \ref{restrict}, \ref{EndAB} and \ref{mapsh} we have the sequence of \Ai-homomorphisms:
$$\left(C,\mu \right)\stackrel{\varphi}{\lto}\left(A\otimes B,\eta\right)\stackrel{F}{\lto}End_{A,B}\stackrel{S}{\longleftarrow}  A\otimes_\infty B.$$
We already saw that $\varphi$ and $S$ are quasi-isomorphisms, so all that remains is to show that $F$ is also a quasi-isomorphism. Recall that the map $\overline P$ introduced in the proof of Proposition \ref{EndAB} satisfies
$$\overline P(F_{1,0}(u\otimes v))=(-1)^{(|u|+|v|)}\eta_{2,0}(u\otimes v, e_A\otimes e_B)= u \otimes v.$$
Hence $\overline P\circ F_{1,0}=id_{A\otimes B}$. Since $\overline P$ is an isomorphism on cohomology, we conclude that $F_{1,0}$ induces an isomorphism in cohomology. This completes the proof of the theorem.
\end{proof}

\section{Bounding cochains on the tensor product}
\subsection{Bounding Cochains}
In this subsection we use the description of the tensor product given in Theorem \ref{filteredtensor} to describe (some of) the bounding cochains on $A\otimes_{\infty} B$. Namely, we will construct a map
$$\boxtimes: MC(A) \times MC(B) \lto MC(A\otimes B, \m^{\otimes}).$$
We begin with the following preliminary lemma:
\begin{lem}\label{trees}
Consider a tree $U \in L_n$, $\beta \in G_A$ and $a_1,\ldots a_n \in A$. Let $j_a$ be the number of $a_i$ such that $a_i=e_A$ and denote by $b(U)$ the number of vertices of $U$ of valency $2$. If $U^A_{\beta}(a_1,\ldots,a_n)\neq 0$, then $b(U)\geq j_a$ unless $j_a=n$ and $\beta=0$, in which case $b(U)=n-1$. 

Moreover the same is true if $U \in R_n$.
\end{lem}
\begin{proof}
Let $w_1^i,\ldots,w^i_{l_i}$ be the set of internal vertices of $P_i$, the path from the $i^{th}$ leaf of $U$ to the root ($w^i_1$ adjacent to the leaf, $w^i_{l_i}$ adjacent to the root). If $a_i=e_A$ then $w^i_1$ must be a binary vertex (of weight zero) since $A$ is unital. This implies $b(U)\geq j_a$ unless there is $1\leq l\leq n$ such that $a_l=e_A$ and $a_{l+1}=e_A$, in which case $\m_{2,0}(a_l,a_{l+1})=e_A$. Hence $w^l_2$ is also binary.

By definition of $L_n$, the edge from $w_1^l$ to $w_2^l$ is right leaning and the other (since $w_2^l$ has valency two) incoming edge at $w_2^l$ is adjacent to a leaf. This implies  $b(U)\geq j_a$ unless $a_{l+2}= e_A$. Iterating this argument we conclude that the only way $b(U)< j_a$ is if $l=1$ and $a_i=e_A$ for all $i$. This implies $U$ is the minimal binary tree, $j_a=n$ and $b(U)=n-1$.

The same argument (interchanging right and left) applied to $U\in R_n$, proves the second claim. The only difference is that $U$ is the maximal binary tree when $j_a=n$.
\end{proof}

\begin{prop}\label{box}
 Let $x \in \widehat{MC}(A)$, $y \in \widehat{MC}(B)$, and define $x\boxtimes y= x\otimes e_B + e_A\otimes y \in \widehat{A\otimes B}_0$. Then $x\boxtimes y \in \widehat{MC}(A\otimes B,\m^\otimes)$ and 
$$\mP(x\boxtimes y)=\mP(x)+\mP(y).$$
\end{prop}
\begin{proof}
From Proposition \ref{alttensor} we have
\begin{align}
\m_{n,\beta_1 \times \beta_2}^\otimes(x\boxtimes y,\ldots,x\boxtimes y)=\sum_{\substack{U, W \\a_i,b_j}}\pm U^A_{\beta_1} (a_1,\ldots,a_n)\otimes W^B_{\beta_2}(b_1,\ldots,b_n),\nonumber
\end{align}
where for each $i$, $a_i=x$ and $b_i=e_B$, or $a_i=e_A$ and $b_i=y$. Moreover $U \in L_n$, $W \in R_n$ and
$\vert E_{int}(U)\vert + \vert E_{int}(W)\vert =n-2$. Equivalently $\vert V(U)\vert +\vert V(W) \vert =n$ and therefore
\begin{align}\label{productedges}
b(U)+b(W) \leq n.
\end{align}
Observe that $j_a + j_b =n$, then if $j_a, j_b \neq 0$,
$$n= j_a +j_b \leq b(U) + b(W) \leq n$$
by the previous lemma. This implies that both $U$ and $W$ are binary and therefore $\vert V(U)\vert =n-1=\vert V(W) \vert$, which contradicts (\ref{productedges}), unless $n=2$. Therefore either $n=2$, or $j_a=n$ and $\beta_1=0$, or $j_b=n$ and $\beta_2=0$. When $n\neq 2$, $U$ and $W$ are (respectively) the minimal and maximal binary trees and we get
\begin{align}
\m^\otimes_{n,\beta_1 \times \beta_2}(x\boxtimes y,\ldots,x\boxtimes y)=\left\{\begin{array}{ll}
\m_{n,\beta_1}(x,\ldots,x)\otimes e_B&,\beta_2=0,\\
e_A\otimes \m_{n,\beta_2}(y,\ldots,y)&,\beta_1=0,\\
0 &, \beta_1,\beta_2\neq0.
\end{array}\right.\nonumber
\end{align} 
When $n=2$, there are two extra terms
$$\m^{\otimes}_{2,\beta_1 \times \beta_2}(x\otimes e_B,e_A \otimes y)+\m^{\otimes}_{2,\beta_1 \times \beta_2}(e_A \otimes y, x\otimes e_B),$$
which cancel because $A$ and $B$ are unital. 
Since $x$ and $y$ are bounding cochains we have 
\begin{align}
\sum_{n} \m^\otimes_{n}(x\boxtimes y,\ldots,x\boxtimes y)&=\big[ \sum_{n} \m_n(x,\ldots,x)\big] \otimes e_B +  e_A \otimes \big[\sum_{n} \m_n(y,\ldots,y)\big]\nonumber\\
 &= \mP(x) e_A \otimes e_B + e_A \otimes \mP(y) e_B \\
 &= (\mP(x)+\mP(y)) e_{A \otimes B}. \nonumber
\end{align}
\end{proof}

Next we will show that the map $\boxtimes$ preserves gauge equivalence. For this we need to make a small digression on models for $A^{[0,1]}$. As in \cite{Sei}, we will construct a model for $A^{[0,1]}$ as the tensor product of $A$ with a filtered dg-algebra $I$. The dg-algebra $I$ is the (normalized) cochain complex of the standard one-simplex.

Consider the graded vector space $I=I^0 \oplus I^1$, with $I^0$ generated by elements $u_0$ and $u_1$ and $I^1$ generated by $h$. We define operations $\mu_{k,\beta}$ on $I$ by setting
\begin{align}
 -\mu_{1,0}(u_0)&= h=\mu_{1,0}(u_1),\nonumber\\
 \mu_{2,0}(u_0,h)=h=-\mu_{2,0}(h,u_1),& \ \  \mu_{2,0}(u_0,u_0)=u_0, \ \ \mu_{2,0}(u_1,u_1)=u_1, \nonumber
\end{align}
and define all other operations to be trivial.

\begin{lem}
 $(I,\mu)$, as defined above, is a $ \{ 0 \}$-gapped filtered dg-algebra
 with unit $e=u_0+u_1$. There are naive quasi-isomorphisms $\bar{i}: \mathbb{K} \lto I$ and $\bar p_0, \bar p_1 : I \lto \mathbb{K}$ satisfying $\bar p_0 \circ \bar{i} = \bar p_1 \circ \bar{i} = id_{\mathbb{K}}$. Additionally the map $\bar p_0 \oplus \bar p_1 : I \lto \mathbb{K} \oplus \mathbb{K}$ is a surjection.
 \end{lem}
\begin{proof}
This is a simple computation. The map $\bar{i}$ is defined as $\bar{i}(1)= u_0 +u_1$ and $\bar p_j$ is defined as the projection to the subspace generated by $u_j$.
\end{proof}

\begin{prop}
 The filtered \Ai-algebra $(A\otimes I, \m^{\otimes})$ together with the naive quasi-isomorphisms $p_j=id \otimes \bar p_j: A\otimes I \lto A\otimes \mathbb{K}=A$ and $i=id \otimes \bar{i}: A=A\otimes \mathbb{K} \lto A\otimes I$ is a model for $A^{[0,1]}$.
\end{prop}
\begin{proof}
We just need to observe that given naive maps $f: A_1 \lto A_2$ and $g: B_1 \lto B_2$, the naive map $f\otimes g : A_1\otimes B_1 \lto A_2 \otimes B_2$ is an \Ai-homomorphism. All the properties are then straightforward.
\end{proof}

\begin{rmk}
 This model for $A^{[0,1]}$ is exactly the one described in \cite[Lemma 4.2.25]{FOOO}. The other explicit model for $A^{[0,1]}$ given in \cite[Lemma 4.2.13]{FOOO} can also be seen to be a tensor product. In fact, $Map([0,1], A)$ is the tensor product of $A$ with the de Rham complex of differential forms on the unit interval $[0,1]$.
\end{rmk}

Using this model for $A^{[0,1]}$ we can describe gauge equivalence explicitly . Two bounding cochains $x_0$ and $x_1$ are gauge equivalent if $\mP(x_0)=\mP(x_1)$ and there is $c=\sum c_i T^{\lambda_i}$ with each $c_i$ of even degree  (or degree zero in the graded case) satisfying
\begin{align}\label{gaugeeq}
 x_0-x_1=\sum_{i,j} \m_{i+j+1}(x_0,\ldots,x_0,c,x_1\ldots,x_1). 
\end{align}

After this digression we are now ready to prove the following proposition:
\begin{prop}
 The map $\boxtimes$ preserves gauge equivalence and descends to a map
 $$\boxtimes: MC(A) \times MC(B) \lto MC(A\otimes_{\infty}B).$$
\end{prop}
\begin{proof}
We need to prove that if $x_0 \sim x_1 \in \widehat{MC}(A)$ and $y_0 \sim y_1 \in \widehat{MC}(B)$ then $x_0\boxtimes y_0 \sim x_1\boxtimes y_1$. We will do this in two steps showing that $x_0\boxtimes y_0 \sim x_1\boxtimes y_0 \sim x_1\boxtimes y_1$.

For the first equivalence we observe that the \Ai-algebra $(A\otimes I)\otimes B$ together with the maps $p'_j=p_j\otimes id: (A\otimes I)\otimes B \lto A \otimes B$ and $i'=i\otimes id: A\otimes B \lto (A\otimes I)\otimes B$ is another model for $(A\otimes B)^{[0,1]}$.
Now by definition there is $c \in \widehat{MC}(A\otimes I)$ such that $p_{j}(c)=x_j$. Then, by Proposition \ref{box}, $c\boxtimes y_0 \in \widehat{MC}((A\otimes I)\otimes B)$ and we have $p'_j(c\boxtimes y_0)= x_j \boxtimes y_0$. Therefore $x_0\boxtimes y_0 \sim x_1\boxtimes y_0$, by definition.

Similarly, using $A\otimes(B \otimes I)$ as a model for $(A\otimes B)^{[0,1]}$, we see that $x_1\boxtimes y_0 \sim x_1\boxtimes y_1$.
\end{proof}

In general, the map $\boxtimes$ is neither injective nor surjective. However, when $A$ and $B$ are graded, there is one simple situation  where we can show that $\boxtimes$ is a bijection.

\begin{defn}
 A filtered \Ai-algebra $(A, \m)$ is connected if $H^0(A, \m_{1,0})= \mathbb{K}$.
\end{defn}

\begin{prop}
 Let $A$ and $B$ be graded, connected, filtered \Ai-algebras. Then $\boxtimes: MC(A) \times MC(B) \lto MC(A\otimes_{\infty}B)$ is a bijection.
\end{prop}
\begin{proof}
First note that by replacing $A$ and $B$ with their canonical models, we can assume that $\m_{1,0}=0$ and thus $A^0= \mathbb{K} e_A$ and $B^0= \mathbb{K} e_B$. Second, we can ignore gauge equivalence since equation (\ref{gaugeeq}) together with the grading and unitality imply that this relation is trivial.

We now proceed by direct computation. From the definition of grading on the tensor product
$$(A\otimes B)^1= e_A \otimes B^1 \oplus A^1\otimes e_B.$$
Thus any element $z \in(A\otimes B)^1 \hat{\otimes} \Lambda_0$ is of the form  $z= x \otimes e_B + e_A \otimes y = x \boxtimes y$. As in the proof of Proposition \ref{box}, we have
$$\sum_{n} \m^\otimes_{n}(z,\ldots,z)=\big[ \sum_{n} \m_n(x,\ldots,x)\big] \otimes e_B +  e_A \otimes \big[\sum_{n} \m_n(y,\ldots,y)\big].$$
Hence $ z \in MC(A \otimes_\infty B)$ if and only if $x \in MC(A)$ and $y \in MC(B)$.
\end{proof}

\subsection{Cohomology and deformations of the tensor product}

As we saw in Section 2, each bounding cochain $ x \in MC(A)$ determines a deformation of the \Ai-algebra $A$. In this section we will describe this deformation for the bounding cochains $x\boxtimes y \in MC(A \otimes_{\infty} B)$ constructed in the previous subsection. 

To avoid working with completions we will assume that all the \Ai-algebras are \emph{compact}. We say an \Ai-algebra is compact if $H^*(A,\m_{1,0})$ is a finite dimensional $\mathbb{K}$-vector space. By replacing the \Ai-algebra by its canonical model, if necessary, this is equivalent to assuming that $A$ is finite dimensional.
Under this assumption, the tensor product $A \otimes_{\mathbb{K}} \Lambda_0$ is already complete, therefore
$$\hat{A}_0=A \otimes_{\mathbb{K}} \Lambda_0 \ \textrm{and} \ \ \hat{A}=A \otimes_{\mathbb{K}} \Lambda.$$
Similarly $\widehat{A\otimes B}_0 = \hat{A}_0 \otimes_{\Lambda_0} \hat{B}_0$ and $\widehat{A\otimes B} = \hat{A} \otimes_{\Lambda} \hat{B}$.

Under these assumptions, we will show that $ (\widehat{A\otimes B}, \m^{\otimes, x\boxtimes y}_k)$ is quasi-isomorphic to $(\hat{A},\m^x_k)\otimes_\infty (\hat{B},\m^y_k)$, as classical \Ai-algebras over $\Lambda$. This will be a direct consequence of Theorem \ref{criterion}, once we prove the following proposition:

\begin{prop}\label{subalgtensor}
 Let $x \in MC(A)$ and $y \in MC(B)$ be bounding cochains and consider $z=x\boxtimes y\in MC(A\otimes_\infty B)$. Then $(\hat{A},\m^x_k)$ and $(\hat{B},\m^y_k)$ are commuting subalgebras of $ (\hat{A} \otimes_{\Lambda} \hat{B}, \m^{\otimes, z}_k)$, under the inclusions
 \begin{align}
  \begin{array}{lll}
\hat{A} \lto \hat{A} \otimes_{\Lambda} \hat{B}, & & \hat{B} \lto \hat{A} \otimes_{\Lambda} \hat{B}\\
a \lto a\otimes e_B                            & &  b \lto e_A\otimes b.\\
\end{array}\nonumber
 \end{align}
\end{prop}
\begin{proof}
First we need to prove that $\hat{A}$ is a subalgebra, that is, for each $n>0$ and $a_1,\ldots,a_n \in \hat{A}$ we need to check
$$\m^{\otimes,z}_n(a_1,\ldots,a_n)= \m^{x}_n(a_1,\ldots,a_n).$$
Using the same argument as in Proposition \ref{box} we obtain
\begin{align}
 \m^{\otimes,z}_n(a_1,\ldots,a_n)&=\sum_{k=i_0+\ldots +i_n} \m_{n+k}(x,\ldots,x,a_1,x,\ldots,x, a_n, x,\ldots,x)\otimes e_B \nonumber\\
 &= \m^{x}_n(a_1,\ldots,a_n).\nonumber
\end{align}

Analogously we see that $B$ is a subalgebra. 

Next we check the commuting relations in Definition \ref{comsubalg}. Consider $c, c_1,\ldots c_n \in \hat{A} \otimes_{\Lambda} \hat{B}$ such that for all $i$, either $c_i=a_i\otimes e_B$ or $c_i=e_A\otimes b_i$ and $c=a\otimes b$. 
Again, repeating the same argument we conclude that
$$\m^{\otimes,z}_n(c_1,\ldots,c_n)=0$$
unless $c_i=a_i\otimes e_B$ for all $i$, $c_i=e_A\otimes b_i$ for all $i$ or  $n=2$. In the latter case, there are two possibilities; we check the case $c_1=a_1\otimes e_B$ and $c_2=e_A \otimes b_2$. We compute
 $$\m^{\otimes,z}_2 (c_1,c_2)=\m_2(a_1, e_A)\otimes \m_2(e_B, b_2) = (-1)^{\vert a_1 \vert} a_1 \otimes b_2,$$
 $$\m^{\otimes,z}_2 (c_2,c_1)=(-1)^{|a_1||b_2|}\m_2(e_A, a_1)\otimes \m_2(b_2, e_B) = (-1)^{\vert b_2 \vert+ |a_1||b_2|} a_1 \otimes b_2.$$
Therefore we conclude $ \m^{\otimes,z}_2 (c_1,c_2) + (-1)^{\vert\vert c_1\vert\vert \vert\vert c_2\vert\vert}\m^{\otimes,z}_2 (c_2,c_1)=0$.

Finally we compute
\begin{align}
 \m^{\otimes,z}_{n+1}(c_1,\ldots ,c_i, c, &c_{i+1}, \ldots , c_n)=\sum_{k=i_0+\ldots +i_n}\m^{\otimes}_{n+k+1}(z,\ldots,z,c_1,z,\ldots,c_n,z,\ldots,z)\nonumber\\
 &=\sum_{k=i_0+\ldots +i_n} \sum _{U, W} \pm U^A \otimes W^B (z,\ldots,z,c_1,z,\ldots,c_n,z,\ldots,z).\nonumber
\end{align}
This sum can be further expanded since $z=x\otimes e_B + e_A \otimes y$. In the notation of Lemma \ref{trees}, we have $k+n = j_{a,x} + j_{b,y}$, therefore the only nontrivial contributions come from trees $U, W$ that satisfy
$$k+n = j_{a,x} + j_{b,y} \leq b(U)+ b(W) \leq n+k+1.$$
Hence there are two possibilities: $U$ is binary and $W$ has a single internal vertex, or vice-versa.
When $n=0$, both cases contribute and we have
\begin{align}
 \m^{\otimes,z}_1(a\otimes b)= \m^x_1(a)\otimes b + (-1)^{\vert a \vert}a \otimes \m^y_1(b).\nonumber
\end{align}
When $n\geq 1$, since $A$ and $B$ are unital, the only contributions come from the first case, when all the $c_i=e_A \otimes b_i$, or the second, when all the $c_i=a_i\otimes e_B$. A simple computation, using the formula for the signs in Theorem \ref{filteredtensor} shows that in the first case we get
\begin{align}
 \m^{\otimes,z}_{n+1}(c_1,\ldots ,c_i, c, \ldots &, c_n)= (-1)^{\vert a \vert (1+\sum_{j \leq i}\vert\vert b_j\vert\vert)} a \otimes \m^{y}_{n+1}(b_1,\ldots,b_i,b,\ldots,b_n)\nonumber
\end{align}
and in the second,
\begin{align}
 \m^{\otimes,z}_{n+1}(c_1,\ldots ,c_i, c, \ldots &, c_n)= (-1)^{\vert b \vert \sum_{i+1 \leq j}\vert\vert a_j\vert\vert} \m^{x}_{n+1}(a_1,\ldots,a_i,a,\ldots,a_n)\otimes b, \nonumber
\end{align}
as required.
\end{proof}

\begin{cor}
 Let $x \in MC(A)$ and $y \in MC(B)$ be bounding cochains and consider the bounding cochain $z=x\boxtimes y\in MC(A\otimes_\infty B)$. We have the isomorphism of classical \Ai-algebras over $\Lambda$,
 $$(\hat{A},\m^x)\otimes_\infty (\hat{B},\m^y) \simeq (\widehat{A\otimes B}, \m^{\otimes, z}).$$
\end{cor}
\begin{proof}
Proposition \ref{subalgtensor} and Theorem \ref{criterion} immediately imply the result, since $K:\hat{A} \otimes_\Lambda \hat{B} \lto \widehat{A\otimes B}$ is simply the identity.
\end{proof}

Proposition \ref{subalgtensor} has one additional consequence. 

\begin{cor} Consider $x \in MC(A)$ and $y \in MC(B)$. For each $n \in \mathbb{Z}_2$, we have the exact sequence
\begin{align}0\lto\bigoplus_{i+j=n}H^i(A, x;\Lambda_0)\otimes H^j(B, y;\Lambda_0)\lto H^n(A \otimes_{\infty} B, x\boxtimes y;\Lambda_0)\nonumber\\
\lto\bigoplus_{i+j=n-1}Tor_1^{\Lambda_0}(H^i(A, x;\Lambda_0),H^j(B, y;\Lambda_0))\lto 0.\nonumber
\end{align}
\end{cor}
\begin{proof}
During the proof of Proposition \ref{subalgtensor} we saw that $\m_1^{\otimes, x\boxtimes y}= \m^x_1 \otimes id + id \otimes \m^y_1$. The result then follows from the usual K\"unneth formula \cite[Theorem 3.6.1]{Wei}, once we prove that $\hat{A}^i$ and $\m_1^{x}(\hat{A}^i)$ are flat $\Lambda_0$-modules, for $i \in \mathbb{Z}_2$. The module $\hat{A}^i$ is free by assumption and $\m_1^{x}(\hat{A}^i)$ is a finitely generated submodule of $\hat{A}^{i+1}$, therefore it is also free by Corollary 2.6.7 in \cite{FOOO}.
\end{proof}

\bibliographystyle{plain}
\bibliography{biblio_tensor}

\medskip

\noindent Address:

\noindent The Mathematical Institute, Radcliffe Observatory Quarter,
Woodstock Road, Oxford, OX2 6GG, U.K.

\noindent E-mail: {\tt camposamorim@maths.ox.ac.uk}

\end{document}